\documentclass[fleqn]{article}
\usepackage{amssymb,latexsym,theorem}
\newtheorem{lemma}{Lemma}[section]
\newtheorem{prop}[lemma]{Proposition}
\newtheorem{theo}[lemma]{Theorem}
\newtheorem{co}[lemma]{Corollary}

\def\pr{\noindent{\bf Proof. }}
\def\eop{\hspace*{\fill}$\Box$}
\title{Unipotent representations of Lie incidence geometries}
\author{Antonio Pasini}
\date{}

\begin{document}
\maketitle

\begin{abstract}
If a geometry $\Gamma$ is isomorphic to the residue of a point $A$ of a shadow geometry of a spherical building $\Delta$, a representation $\varepsilon_\Delta^A$ of $\Gamma$ can be given in the unipotent radical $U_{A^*}$ of the stabilizer in $\mathrm{Aut}(\Delta)$ of a flag $A^*$ of $\Delta$ opposite to $A$, every element of $\Gamma$ being mapped onto a suitable subgroup of $U_{A^*}$. We call such a representation a unipotent representation. We develope some theory for unipotent representations and we examine a number of interesting cases, where a projective embedding of a Lie incidence geometry $\Gamma$ can be obtained as a quotient of a suitable unipotent representation $\varepsilon_\Delta^A$ by factorizing over the derived subgroup of $U_{A^*}$, while $\varepsilon^A_\Delta$ itself is not a proper quotient of any other representation of $\Gamma$.     
\end{abstract}

\section{Introduction}

A rich literature exists on projective embeddings of point-line geometries, but non-projective embeddings have been considered too, where a geometry is embedded in a group rather than a vector space. For instance, various interesting families of geometries exist with lines of size $3$, associated to sporadic simple groups. A full understanding of these geometries and the groups acting on them requires a thorough investigation of their representations (Ivanov and Shpectorov \cite{IS}). In the terminology of this paper, a representation as defined in \cite{IS} is just a locally projective but possibly non-projective embedding. Still in the investigation of geometries for finite simple groups, geometries where all lines have $p+1$ points for a small prime $p > 2$ also occur. Representations can be considered for them too (Ivanov \cite{Sasha2001}). 

Non-projective embeddings can also be interesting in perspectives different from the  above. For instance, starting from a projective embedding $\varepsilon$ of an incidence geometry $\Gamma$ of Lie type we can look for the largest possibly non-projective embedding $\tilde{\varepsilon}$ admitting $\varepsilon$ as a quotient. We call $\tilde{\varepsilon}$ the hull of $\varepsilon$. Knowing the hull $\tilde{\varepsilon}$ of $\varepsilon$ is the same as knowing the universal cover of a geometry associated to $\varepsilon$, which we call the expansion of $\Gamma$ by $\varepsilon$. In many interesting cases that universal cover is the geometry far from an element of a suitable spherical building $\Delta$ containing $\Gamma$ as a residue. So, non-projective embeddings can also be useful to better understand certain aspects of spherical buildings.  

We have moved some steps in this trend in \cite{PasEE}. In the present paper we shall continue that investigation. In particular, let $\Delta$ be a spherical building such that $\Gamma$ is isomorphic to the residue of a point $A$ of the $J$-shadow geometry $\Delta_J$ of $\Delta$, for a subset $J$ of the set of types of $\Delta$. We define an embedding $\varepsilon_\Delta^J$ of $\Gamma$ in the unipotent radical $U_{A^*}$ of the stabilizer of $A^*$ in $\mathrm{Aut}(\Delta)$, for a flag $A^*$ of $\Delta$ opposite to $A$. We call $\varepsilon^J_\Delta$ a unipotent representation of $\Gamma$. In many cases, given a projective embedding $\varepsilon$ of $\Gamma$, defined independently of $\Delta$, we can choose the pair $(\Delta,J)$ in such a way that $\varepsilon^J_\Delta$ is the hull of $\varepsilon$. In short, $\tilde{\varepsilon}\cong \varepsilon_\Delta^J$. 

A number of situations like this are considered in \cite[Sections 6-9]{PasEE}, but in that paper unipotent representations were not put in evidence as they deserved. In the present paper we put them in the right light. We develope some general theory for unipotent representations, we revisit a number of results obtained in \cite[Sections 6-9]{PasEE} and examine two cases that have not been considered in \cite{PasEE}, where $\Gamma = E_{6,1}(\mathbb{F})$ or $\Gamma = E_{7,7}(p)$ for a prime $p$ (notation as in Section 9 of this paper) and $\varepsilon$ is the $27$-dimensional projective embedding of $E_{6,1}(\mathbb{F})$ or the $56$-dimensional projective embedding of $E_{7,7}(p)$ respctively. We prove that the isomorphism $\tilde{\varepsilon}\cong \varepsilon_\Delta^J$ also holds in these two cases, with $(\Delta,J) = (E_{7}(\mathbb{F}),\{7\})$ when $\Gamma = E_{6,1}(\mathbb{F})$ and $(\Delta,J) = (E_{8}(p),\{8\})$ for $\Gamma = E_{7,7}(p)$. As a by-product of this result we obtain that the subgeometry of $E_7(\mathbb{F})$ far from a $7$-element and the subgeometry of $E_8(p)$ far from an $8$-element are simply connected.\\

\noindent
{\bf Organization of the paper.} Basics on embeddings, theirs hulls and expansions, shadow geometries and subgeometries far from a flag are recalled in Section 2. In Section 3 we define unipotent representations and develope some theory for them. In the last subsection of Section 3 we examine various unipotent representations of $\mathrm{PG}(n-1,\mathbb{K})$. In particular, we show that quadratic veronesean embeddings and hermitian veronesean embeddings of $\mathrm{PG}(n-1,\mathbb{K})$ are in fact unipotent representations.

A number of results on hulls of projective embeddings of particular Lie incidence geometries are collected in Section 4. The hull of each of the embeddings considered in Section 4 is a unipotent representation. All theorems of Section 4 but two have been proved in \cite[Sections 7-9]{PasEE}. As said above, the two new theorems deal with with $E_{6,1}(\mathbb{F})$ and $E_{7,7}(p)$. They are proved in Sections 7 and 8 respectively. Section 6 contains some generalities on embeddings and subspaces of geometries, to be exploited in Sections 7 and 8. A few conjectures and problems suggested by the results of Section 4 are proposed in Section 5. 

The notation we use to denote particular buildings of sperical type and their shadow geometries is quite standard. Anyway,
we recall it in Section 9. \\

\noindent
{\bf Warning.} Basic notions from diagram geometry will be freely used throughout this paper. We presume that the reader is familiar with them. If not, we refer to Pasini \cite{PasDG} or Buekenhout and Cohen \cite[Chapters 1-3]{BC}. 

\section{Embeddings, shadows and opposition} 

In this section we recall some basics on embeddings of point-line geometries in groups and on shadows and opposition in buildings, and we fix some terminology and notation to be used in the rest of the paper.

\subsection{Embeddings} 

Non-projective embeddings can be defined for any geometry of rank $n$ with a string-shaped diagram, as in \cite{PasEE} (see also Ivanov \cite{Sahsa2001}), but in this paper we prefer to stick to point-line geometries. In this paper a {\em point-line geometry} is a pair $\Gamma = ({\cal P},{\cal L})$ where $\cal P$ (the {\em set of points}) is a non-empty set and $\cal L$ is a family of subsets of $\cal P$, called {\em lines}, satisfying all the following: every line contains at least two points, every point belongs to at least two lines, no two lines have more than one point in common and the collinearity graph of $\Gamma$ is connected. Actually, this definition of point-line geometry is more restrictive than other definitions commonly used in the literature, but it fits our needs here. 

The points and the lines of $\Gamma$ will be called {\em elements} of $\Gamma$. Accordingly, we write $X\in\Gamma$ for short instead of $X\in{\cal P}\cup{\cal L}$.    

An {\em embedding} of a point-line geometry $\Gamma = ({\cal P},{\cal L})$ into a group $G$ is a mapping from ${\cal P}\cup{\cal L}$ into the subgroup lattice of $G$ satisfying all the following:

\begin{itemize}
\item[$(\mathrm{E}1)$] $\varepsilon(p)\cap\varepsilon(q) = 1 < \varepsilon(p)$ for any two distinct points $p, q\in {\cal P}$.  
\item[$(\mathrm{E}2)$] For a point $p$ and a line $L$, we have $\varepsilon(p)\leq \varepsilon(L)$ if and only if $p\in L$.
\item[$(\mathrm{E}3)$] $\varepsilon(L) = \langle \varepsilon(p)\rangle_{p\in L}$ for every line $L\in{\cal L}$. 
\item[$(\mathrm{E}4)$] $G = \langle \varepsilon(p)\rangle_{p\in{\cal P}}$.
\end{itemize}
Condition $(\mathrm{E}1)$ is slightly stronger than the correspondent condition implicit in the definition given in \cite{PasEE}, but it holds in all cases we shall consider in this paper. Note that $(\mathrm{E}1)$ and $(\mathrm{E}2)$ imply that, for any two distinct lines $L, M\in{\cal L}$, the intersection $\varepsilon(L)\cap\varepsilon(M)$ is properly contained in either of $\varepsilon(L)$ and $\varepsilon(M)$. By this fact and $(\mathrm{E}1)$, the mapping $\varepsilon$ is injective. 
 
Henceforth when writing $\varepsilon:\Gamma\rightarrow G$ we mean that $\varepsilon$ is an embedding of $\Gamma$ into $G$. We call $G$ the {\em codomain} of $\varepsilon$. 

If $G$ is abelian then $\varepsilon$ is said to be {\em abelian}. Suppose that $G$ is the additive group of a vector space $V$ and that $\varepsilon(X)$ is (the additive group of) a vector subspace of $V$, for every $X\in \Gamma$. Then we say that the embedding $\varepsilon$ is {\em linear} and {\em defined over} $\mathbb{K}$, where $\mathbb{K}$ is the underlying division ring of $V$. By a slight modification of our previous conventions, we take $V$ as the {\em codomain} of $\varepsilon$, thus writing $\varepsilon:\Gamma\rightarrow V$ instead of $\varepsilon:\Gamma\rightarrow G$. The dimension $\mathrm{dim}(V)$ is called the {\em dimension} of $\varepsilon$ and denoted by $\mathrm{dim}(\varepsilon)$. 

A linear embedding $\varepsilon:\Gamma\rightarrow V$ is said to be {\em projective} if the following holds, where for an element $X\in\Gamma$ the symbol $V_X$ stands for $\varepsilon(X)$, but regarded as a subspace of $V$. 

\begin{itemize} 
\item[$(\mathrm{PE})$] For every line $L$, the vector space $V_L$ is $2$-dimensional and $\{V_p\}_{p\in L}$ is the collection of all $1$-dimensional subspaces of $V_L$.
\end{itemize}
Note that, since every point of $\Gamma$ belongs to at least one line (in fact at least two), $(\mathrm{PE})$ implies that $V_p$ is $1$-dimensional, for every point $p$. 

If $\Gamma$ admits a projective embedding and all projective embeddings of $\Gamma$ are defined over the same division ring $\mathbb{K}$ (as when $\Gamma$ is finite or contains a subspace isomorphic to a projective plane or a thick generalized quadrangle, for instance) then $\Gamma$ is said to be {\em defined over} $\mathbb{K}$. 

The previous definitions can be weakened as follows. Without assuming anything on $G$,  suppose that for a given division ring $\mathbb{K}$ and every element $X\in\Gamma$ a $\mathbb{K}$-vector space $V_X$ exists such that $\varepsilon(X)$ is the additive group of $V_X$ and, if a point $p$ bellongs to a line $L$ then $V_p$ is a vector subspace of $V_L$. Then we say that $\varepsilon$ is {\em locally linear} and {\em defined over} $\mathbb{K}$. If moreover $(\mathrm{PE})$ holds, then $\varepsilon$ is said to be {\em locally projective}. 

If all lines of $\Gamma$ have $p+1$ points for a prime $p$ then the locally projective embeddings of $\Gamma$ are just the representations of $\Gamma$ as defined in Ivanov and Shpectorov \cite{IS} (for $p = 2$) and Ivanov \cite{Sasha2001} (for any $p$).  

\subsubsection{Morphisms and quotients}\label{morphisms quotients}

Given two embeddings $\varepsilon_1:\Gamma\rightarrow G_1$ and $\varepsilon_2:\Gamma\rightarrow G_2$, a {\em morphism} $f:\varepsilon_1\rightarrow \varepsilon_2$ from $\varepsilon_1$ to $\varepsilon_2$ is a group-homomorphism $f:G_1\rightarrow G_2$ such that for every element $X\in\Gamma$ the restriction of $f$ to $\varepsilon_1(X)$ is an isomorphism from $\varepsilon_1(X)$ to $\varepsilon_2(X)$. Note that $(\mathrm{E}4)$ forces $f$ to be surjective. If $f$ is injective (whence it is an isomorphism from $G_1$ to $G_2$) then we say that $f$ is an {\em isomorphism} from $\varepsilon_1$ to $\varepsilon_2$. In this case we say that $\varepsilon_1$ and $\varepsilon_2$ are {\em isomorphic} and we write $\varepsilon_1\cong\varepsilon_2$. 

Let $\varepsilon:\Gamma\rightarrow G$ be an embedding and let $U$ be a normal subgroup of $G$ satisfying the following:

\begin{itemize}
\item[$(\mathrm{Q}1)$] $U\cap \varepsilon(X) = 1$ for every element $X\in\Gamma$. 
\item[$(\mathrm{Q}2)$] $(\varepsilon(p)U)\cap(\varepsilon(q)U) = U$ for any two non-collinear points $p,  q\in {\cal P}$.
\end{itemize}
(Note that $(\mathrm{Q}1)$ implies that $(\varepsilon(p)U)\cap(\varepsilon(q)) = U$ for any two distinct collinear points $p$ and $q$.) Then an embedding $\varepsilon/U:\Gamma\rightarrow G/U$ can be defined by setting $(\varepsilon/U)(X) = (\varepsilon(X)U)/U$ for every element $X\in \Gamma$. The canonical projection $\pi_U:G\rightarrow G/U$ is a morphism from $\varepsilon$ to $\varepsilon/U$. We say that $U$ {\em defines a quotient} of $\varepsilon$ and we call $\varepsilon/U$ the {\em quotient} of $\varepsilon$ by $U$.

Clearly, if $\varepsilon_1:\Gamma\rightarrow G_1$ and $\varepsilon_2:\Gamma\rightarrow G_2$ are embeddings and $f:\varepsilon_1\rightarrow\varepsilon_2$ is a morphism then $\mathrm{Ker}(f)$ defines a quotient of $\varepsilon_1$ and $\varepsilon_1/\mathrm{Ker}(f)\cong \varepsilon_2$. In view of this, we take the liberty to say that $\varepsilon_2$ is a {\em quotient} of $\varepsilon_1$.

\subsubsection{Hulls and dominant embeddings}\label{Hull}

Given a point-line geometry $\Gamma = ({\cal P},{\cal L})$ and an embedding $\varepsilon:\Gamma\rightarrow G$, let
\[{\cal A}_\varepsilon := (\{\varepsilon(X)\}_{X\in \Gamma}, \{\iota_{p,L}\}_{p\in{\cal P}, L\in{\cal L}, p\in L})\]
where for a point-line flag $\{p,L\}$ of $\Gamma$ we denote by $\iota_{p,L}$ the inclusion mapping of $\varepsilon(p)$ in $\varepsilon(L)$. Regarding the groups $\varepsilon(X)$ as abstract groups rather than subgroups of $G$, let $\bf U$ be the disjoint union of the sets $\varepsilon(X)\setminus\{1\}$ for $X\in\Gamma$ and $F({\bf U})$ the free group over $\bf U$. Let $W$ be the minimal normal subgroup of $F(\bf U)$ containing all words $x^{-1}\iota_{p,L}(x)$ for every point-line flag $\{p,L\}$ and $x\in \varepsilon(p)\setminus\{1\}$ and all words $x_1^{k_1}...x_n^{k_n}$ with $x_1,..., x_n\in \varepsilon(X)\setminus\{1\}$ and $x_1^{k_1}...x_n^{k_n} = 1$ in $\varepsilon(X)$ for $X\in\Gamma$. Let $U(\varepsilon) := F({\bf U})/W$. In short, $U(\varepsilon)$ is the universal completion of the amalgam ${\cal A}_\varepsilon$.

The function $\tilde{e}$ mapping every element $X\in \Gamma$ onto the subgroup $\varepsilon(X)$ of $U(\varepsilon)$ is an embedding of $\Gamma$ in $U(\varepsilon)$ and the canonical projection $\pi_\varepsilon$ of $U(\varepsilon)$ onto $G$ is a morphism from $\tilde{\varepsilon}$ onto $\varepsilon$. We call $\tilde{\varepsilon}$ the {\em hull} of $\varepsilon$. It satisfies the following `universal' property:

\begin{itemize}
\item[$(\mathrm{U})$] For any embedding $\varepsilon'$ of $\Gamma$, if $\varepsilon$ is a quotient of $\varepsilon'$ then for every morphism $f:\varepsilon'\rightarrow \varepsilon$ a unique morphism $g:\tilde{\varepsilon}\rightarrow \varepsilon'$ exists such that $f\cdot g = \pi_\varepsilon$. 
\end{itemize}
Property $(\mathrm{U})$ uniquely determines $\tilde{\varepsilon}$ up to isomorphisms. Borrowing a word from Tits \cite[8.5.2]{Tits}, we say that the embedding  $\varepsilon$ is {\em dominant} if $\pi_\varepsilon$ is an isomorphism. Equivalently, $\varepsilon\cong \tilde{\varepsilon}$. 
 
Suppose that $\varepsilon$ is abelian. Then the commutator subgroup $U(\varepsilon)'$ of $U(\varepsilon)$ defines a quotient $\tilde{\varepsilon}_{\mathrm{ab}} := \tilde{\varepsilon}/U(\varepsilon)'$ of $\tilde{\varepsilon}$. We call $\tilde{\varepsilon}_{\mathrm{ab}}$ the {\em abelian hull} of $\varepsilon$.  We also put $U(\varepsilon)_{\mathrm{ab}} := U(\varepsilon)/U(\varepsilon)'$, for short. 

The abelian hull $\tilde{\varepsilon}_{\mathrm{ab}}$ is characterized by a universal property quite similar to $(\mathrm{U})$, except that only abelian embeddings are considered in it and $\pi_\varepsilon$ is replaced with the projection of $U(\varepsilon)_{\mathrm{ab}}$ onto $G$ induced by $\pi_\varepsilon$. 

Suppose moreover that $\varepsilon$ is $\mathbb{K}$-linear for a given division ring $\mathbb{K}$. Hence $G$ is the additive group of a $\mathbb{K}$-vector space $V$ and $\varepsilon(X)$ is a subspace of $V$ for every $X\in\Gamma$.  Let $V$ be $W$ be the subgroup of $U(\varepsilon)_{\mathrm{ab}}$ generated by all sums $\sum_{i=1}^nta_i$, for every scalar $t\in\mathbb{K}$ and every choice of elements $a_1,..., a_n\in\cup_{X\in\Gamma}\varepsilon(X)$ such that $\sum_{i=1}^na_i = 0$ in $U(\varepsilon)_{\mathrm{ab}}$. Needless to say, we are adopting the additive notation for the grouyp $U(\varepsilon)_{\mathrm{ab}}$ and for $i = 1,..., n$ the product $ta_i$ is computed in the vector space $V_X$ supported by $\varepsilon(X)$, for $X\in \Gamma$ such that $a_i\in\varepsilon(X)$. With $W$ defined in this way,  the factor group $U(\varepsilon)_{\mathrm{lin}} := U(\varepsilon)_{\mathrm{ab}}/W$ gets the structure of a $\mathbb{K}$-vector space and the projection $\pi_{\varepsilon,\mathrm{lin}}:U(\varepsilon)_{\mathrm{lin}}\rightarrow G$ induced by $\pi_\varepsilon$ is a linear mapping from the $\mathbb{K}$-vector space $U(\varepsilon)_{\mathrm{lin}}$ to $V$. The subgroup $W$ defines a quotient $\tilde{\varepsilon}_{\mathrm{lin}} := \tilde{\varepsilon}_{\mathrm{ab}}/W$ of $\tilde{\varepsilon}_{\mathrm{ab}}$. We call $\tilde{\varepsilon}_{\mathrm{lin}}$ the {\em linear hull} of $\varepsilon$. The linear hull $\tilde{\varepsilon}_{\mathrm{lin}}$ is characterized by a universal property similar to $(\mathrm{U})$, except that only linear embeddings are considered in itan only semi-linear mappings are taken as morphisms of embeddings. In particular, the role of $\pi_\varepsilon$ is taken over by $\pi_{\varepsilon,\mathrm{lin}}$. We say that a linear embedding is {\em linearly dominant} if it is its own linear hull.

Clearly, the hull of a locally linear (in particular, locally projective) embedding is still locally linear (locally projective). Hence if $\varepsilon$ is projective embedding then its linear hull is projective as well. It is just the projective embedding universal relatively to $\varepsilon$ as defined in Ronan \cite{Ronan} and Shult \cite{Shult} (compare our construction of $\tilde{\varepsilon}_{\mathrm{lin}}$ with the construction by Ronan \cite[Section 2]{Ronan}). Consequently, a projective embedding is linearly dominant if and only if it is relatively universal in the sense of Shult \cite{Shult}.

If a projective embedding is dominant then it is linearly dominant, but the converse is false in generale: the hull of a linearly dominant projective embedding is seldom abelian, let alone linear.

\subsubsection{Absolute projective embeddings}\label{absolute}

Following Kasikova and Shult \cite{KS}, we say that a projective emedding $\varepsilon:\Gamma\rightarrow V$ is {\em absolute} if it is the linear hull of all projective embeddings of $\Gamma$ defined over the same division ring as $\varepsilon$. Linear hulls are uniquely determined modulo isomorphisms. Hence the absolute projective embedding (defined over a given division ring), if it exists, is unique up to isomorphisms. Not every geometry that admits a projective embedding also admits the absolute projective embedding, but many of them do. In particular, each of the Lie incidence geometries to be considered in Section 4 (is defined over a suitable division ring and) admits the absolute projective embedding.  

\subsubsection{Expansions}

Let $\varepsilon:\Gamma\rightarrow G$ be an embedding. The {\em expansion} of $\Gamma$ to $G$ by $\varepsilon$ is the rank $3$ geometry $\mathrm{Ex}(\varepsilon)$ defined as follows. The triple $\{0,1,2\}$ is taken as the set of types of $\mathrm{Ex}(\varepsilon)$. The elements of $\mathrm{Ex}(\varepsilon)$ of type $0$ (called {\em points} of $\mathrm{Ex}(\varepsilon)$) are the elements of $G$. The elements of type $1$ ({\em lines}) and those of type $2$ ({\em planes}) are the cosets $g\cdot\varepsilon(X)$, for $X$ a point or a line of $\Gamma$ respectively. Set theoretic inclusion is taken as the incidence relation. 

The structure $\mathrm{Ex}(\varepsilon)$ defined as above is indeed a residually connected geometry. Its $\{1,2\}$-residues are isomorphic to $\Gamma$ while its $\{0,1\}$ residues are (possibly infinite) nets (affine planes when $\varepsilon$ is projective). The group $G$ in its action on itself by left multiplication yields a group of automorphisms of $\mathrm{Ex}(\varepsilon)$, sharply transitive on the set of points of $\mathrm{Ex}(\varepsilon)$.

If $\bar{\varepsilon}:\Gamma\rightarrow \overline{G}$ is another embedding of $\Gamma$ and $f:\overline{G}\rightarrow G$ a morphism from $\bar{\varepsilon}$ to $\varepsilon$, then $f$ naturally induces 
a morphism $\mathrm{Ex}(f)$ from $\mathrm{Ex}(\bar{\varepsilon})$ to $\mathrm{Ex}(\varepsilon)$, mapping $g\in \overline{G}$ onto $f(g)$ and $g\cdot\bar{\varepsilon}(X)$ onto $f(g)\cdot\varepsilon(X)$ for every $X\in\Gamma$. It is easy to see that $\mathrm{Ex}(f)$ is a covering.  

Let $\widetilde{\mathrm{Ex}(\varepsilon)}$ be the universal cover of $\mathrm{Ex}(\varepsilon)$ and $\tilde{\varepsilon}:\Gamma\rightarrow U(\varepsilon)$ the hull of $\varepsilon$. Let $\pi_\varepsilon:U(\varepsilon)\rightarrow G$ be the projection of $\tilde{\varepsilon}$ onto $\varepsilon$. The following is proved in \cite{PasEE} (see also Ivanov \cite[Lemma 3.9]{Sasha2001} for the special case where all lines of $\Gamma$ have size $p+1$ for a given prime $p$ and $\varepsilon$ is locally projective). 

\begin{prop}\label{expansion}
We have $\widetilde{\mathrm{Ex}(\varepsilon)} = \mathrm{Ex}(\tilde{\varepsilon})$. Moreover, $\mathrm{Ker}(\pi_\varepsilon)$ is the group of deck transformations of the covering $\mathrm{Ex}(\pi_\varepsilon)$, hence it is isomorphic to the homotopy group of $\mathrm{Ex}(\varepsilon)$. 
\end{prop}

\subsubsection{Homogeneity}\label{Homogeneity}

Given an embedding $\varepsilon:\Gamma\rightarrow G$, we say that an automorphism $\alpha$ of $\Gamma$ {\em lifts} to $G$ through $\varepsilon$ if an automorphism $\hat{\alpha}$ of $G$ exists, uniquely determined modulo automorphisms of $\varepsilon$, such that $\hat{\alpha}\cdot \varepsilon = \varepsilon\cdot \alpha$. If all elements of a subgroup $A\leq \mathrm{Aut}(\Gamma)$ lift to $G$ then we say that $A$ {\em lifts} to $G$ through $\varepsilon$ and that $\varepsilon$ is $A$-{\em homogeneous}. The set $\widehat{A}$ of all liftings $\hat{\alpha}$ for $\alpha\in A$ is a subgroup of $\mathrm{Aut}(G)$, called the {\em lifting} of $A$ to $G$. If $\varepsilon$ is $Aut(\Gamma)$-homogeneous then we say that is {\em homogeneous}, for short.  

Let $\varepsilon$ be $A$-homogeneous for a subgroup $A\leq \mathrm{Aut}(\Gamma)$. Then its hull is also $A$-homogenous. Let $U\unlhd G$ define a quotient of $\varepsilon$ and suppose that every $\alpha\in A$ admits a lifting $\hat{\alpha}$ stabilizing $U$. Then $\varepsilon/U$ is $A$-homogeneous. 

\subsection{Shadows and opposition} 

Throughout this subsection $\Delta$ is a thick building of spherical type and rank $n \geq 2$, $I$ is the set of types of $\Delta$ and $\tau:\Delta\rightarrow I$ is its type-function. 

\subsubsection{Shadows}\label{shadows} 

We recall the definition of shadow geometries, as stated in Tits \cite[chapter 12]{Tits} (see also \cite[Chapter 5]{PasDG}, where shadow geometries are called Grassmann geometries and are defined in a more general setting, for any geometry). 

Let $J$ be a subset of $I$ containing at least one type of every irreducible component of the Coxeter diagram of $\Delta$. The $J$-{\em shadow} $\mathrm{sh}_J(F)$ of a nonempty flag $F$ of $\Delta$ is the set of flags of type $J$ incident to $F$. In general, different flags can have the same $J$-shadow. However, given a $J$-shadow $X$, the family ${\cal F}_X$ of flags $F$ such that $\mathrm{sh}_J(F) = X$ admits a smallest member (Tits \cite{Tits}). We denote it by $F_X$. 

Let ${\cal D}$ be the diagram graph of $\Delta$ (see \cite[Chapter 4]{PasDG}, also Buekenhout and Cohen \cite[2.1]{BC}, where diagrams graphs are called digon diagrams). In short, $\cal D$ is the graph obtained from the Coxeter diagram of $\Delta$ by replacing multiple strokes with simple strokes. With $X$ as above and $F\in{\cal F}_X$, let ${\cal D}_F$ be the graph induced by ${\cal D}$ on the cotype $I\setminus\tau(F)$ of $F$ and let ${\cal D}_{F}^J$ be the union of the connected components of ${\cal D}_F$ meeting $J$ non-trivially. (We warn that ${\cal D}_F^J$ might be empty, but it is empty if and only if $X$ is a single $J$-flag, namely $\tau(F)\supseteq J$.) We have ${\cal D}_F^J = {\cal D}_{F_X}^J$ for any $F\in {\cal F}_X$ (Tits \cite{Tits}). This fact allows us to slightly change our notation, writing ${\cal D}_X$ for ${\cal D}_F^J$. The type $\tau(F_X)$ of $F_X$ is the neigborhood of ${\cal D}_X$ in $\cal D$. Clearly, for two $J$-shadows $X$ and $Y$ we have ${\cal D}_X\subseteq {\cal  D}_Y$ if and only if $\tau(F_X)$ separates $\tau(F_Y)$ from $J$, namely every path of $\cal D$ from $J$ to $\tau(F_Y)$ meets $\tau(F_X)$ at some vertex. We have $X\subseteq Y$ if and only if ${\cal D}_X\subseteq{\cal D}_Y$ and $F_X\cup F_Y$ is a flag (Tits \cite{Tits}).  

We take the integer $\tau_J(X) := |{\cal D}_X|+1$ as the {\em type} of the $J$-shadow $X$. Moreover, we say that two $J$-shadows are {\em incident} if one of them contains the other one. In this way the set of $J$-shadows of nonempty flags of $\Delta$ is turned into a residually connected geometry of rank $n$ over the set of types $\{1, 2,,...,n\}$, henceforth denoted by $\Delta_J$ and called the $J$-{\em shadow geometry} of $\Delta$. 

In general, $\Delta_J$ does not belong to any Coxeter diagram, although all residues of $\Delta_J$ of rank 2 are generalized polygons. Indeed it can happen that for some $i < n$ not all $\{i,i+1\}$-residues of $\Delta_J$ have the same gonality. Nevertheless $\Delta_J$ admits a Buekenhout diagram, which has the shape of a string (namely the diagram graph of $\Delta_J$ is a string), the type $1$ corresponding to one of the two end nodes of the string. To fix ideas, let $1$ correspond to the leftmost node. Then the remaining types label the nodes of the diagram in increasing order, from left to right. 

The elements of $\Delta_J$ of type $1$, $2$ and $3$ (if $n > 2$) are called {\em points}, {\em lines} and {\em planes} respectively. Clearly, the points of $\Delta_J$ are just the flags of type $J$. 

We have defined $\Delta_J$ as a geometry of rank $n$, but when $n > 2$ we are often only interested in its point-line system, namely its $\{1,2\}$-truncation. For the sake of pedantry, we should introduce different symbols for $\Delta_J$ and its $\{1,2\}$-truncation, but we prefer to use the symbol $\Delta_J$ for either of them. In any case the context or suitable warnings will make it clear if we regard $\Delta_J$ as a geometry of rank $n$ or a point-line geometry.

In many cases considered in the literature $J$ is a singleton, say $J = \{k\}$. In this case we write $\Delta_k$ for short instead of $\Delta_{\{k\}}$.   

\subsubsection{Local geometries of $\Delta_J$}\label{local}

Suppose now that $n > 2$. Given a point $A$ of $\Delta_J$, namely a flag $A$ of type $J$, let $\Delta_J(A)$ be the $\{2,3\}$-truncation of the residue $\mathrm{Res}_{\Delta_J}(A)$ of $A$ in $\Delta_J$, namely the point-line geometry the points and the lines of which are the lines and the planes of $\Delta_J$ incident to $A$. We denote the point-set and the line-set of $\Delta_J(A)$ by ${\cal P}_A$ and ${\cal L}_A$ respectively and we call $\Delta_J(A)$ the {\em local geometry} of $\Delta_J$ at $A$. 

Given a line $X\in{\cal L}_A$ of $\Delta_J(A)$, we denote by $P(X)$ the set of points of $\Delta_J(A)$ incident to $X$. It is not difficult to see that if $X$ and $Y$ are distinct points of $\Delta_J(A)$ then $X\cap Y = \{A\}$. If $X$ and $Y$ are distinct lines of $\Delta_J(A)$ then either $P(X)\cap P(Y) = \emptyset$ and $X\cap Y = \{A\}$ or $P(X)\cap P(Y) = \{Z\}$ for a unique point $Z\in {\cal P}_A$ and $X\cap Y = Z$. If $X$ is a point and $Y$ is a line then either $X\in P(Y)$ (namely $X\subset Y$) or $X\cap Y = \{A\}$. 

As $n > 2$, the building $\Delta$ is flag-transitive (Tits \cite[Chapter 3]{Tits}). Hence the isomorphism type of $\Delta_J(A)$ does not depend on the choice of the $J$-flag $A$. Note also that when $J$ is a singleton, say $J = \{k\}$, the geometry $\Delta_J(A)$ is isomorphic to the point-line system of the $K$-shadow geometry $\mathrm{Res}_\Delta(A)_K$ of the residue $\mathrm{Res}_\Delta(A)$ of $A$ in $\Delta$, where $K$ is the neighborhood of $k$ in $\cal D$.

\subsubsection{Lower residues of $\Delta_J$}\label{shadows residues sec} 

Given a $J$-shadow $X\in \Delta_J$ with $\tau_J(X)\geq 2$, let $F\in {\cal F}_X$ (notation as in Subsection \ref{shadows}) and let $\mathrm{Res}^J_\Delta(F)$ be the subgeometry of $\mathrm{Res}_\Delta(F)$ formed by the elements of type $i\in{\cal D}_X$. Note that $\mathrm{Res}^J_\Delta(F)$ does not depend on the particular choice of $F$ in ${\cal F}_X$. In particular, if $F_X$ is the minimal element of ${\cal F}_X$ and $\overline{F}\in {\cal F}_X$ is maximal (namely $\tau(\overline{F}) = I\setminus{\cal D}_X$), then $\mathrm{Res}^J_\Delta(F_X) = \mathrm{Res}_\Delta^J(\overline{F}) = \mathrm{Res}_\Delta(\overline{F})$. So, $\mathrm{Res}_\Delta^J(F)$ is a ${\cal D}_X$-residue of $\Delta$. Hence it is a building of spherical type. Obviously, $J\setminus\tau(F)$ meets every conponent of ${\cal D}_X$ non-trivially. Let $\mathrm{Res}_{\Delta_J}^-(X)$ be the subgeometry of $\mathrm{Res}_{\Delta_J}(X)$ formed by the elements of $\mathrm{Res}_{\Delta_J}(X)$ of type less than $\tau_J(X)$. The following is straightforward.

\begin{prop}\label{shadows-residues}
The geometry $\mathrm{Res}_{\Delta_J}^-(X)$ is isomorphic with the $(J\setminus\tau(F))$-shadow geometry 
$\mathrm{Res}^J_\Delta(F)_{J\setminus\tau(F)}$ of the building $\mathrm{Res}^J_\Delta(F)$. 
\end{prop}
  
\subsubsection{Projections}\label{proj}

For a flag $F$ of $\Delta$, let ${\cal C}(F)$ be the set of chambers of $\Delta$ that contain $F$. It is well known that for every chamber $C$ of $\Delta$ there exist a unique chamber $C'\in{\cal C}(F)$ at minimal distance from $C$ (Tits \cite[2.30]{Tits}), distances between chambers being computed in the chamber system of $\Delta$. Moreover $d(C,D) = d(C,C') + d(C',D)$ for any $D\in{\cal C}(F)$. The chamber $C'$ is called the {\em gate} of ${\cal C}(F)$ to $C$. 

Given another flag $A$ of $\Delta$ let ${\cal C}_A(F)$ be the set of gates of ${\cal C}(F)$ to chambers of ${\cal C}(A)$ and put $A_F := \cap_{C\in{\cal C}_A(F)}C$. Clearly, $A_F$ is a flag containing $F$. So, $\mathrm{pr}_F(A) := A_F\setminus F$ is a (possibly empty) flag of $\mathrm{Res}_\Delta(F)$. We call $\mathrm{pr}_F(A)$ the {\em projection} of $A$ onto $\mathrm{Res}_\Delta(F)$.  

\subsubsection{Opposition and geometries far from a flag}\label{Far}   

Two chambers of the building $\Delta$ are said to be {\em opposite} if they have maximal distance in the chamber system of $\Delta$. Two flags $F$ and $F^*$ are {\em opposite} if every chamber containing $F$ or $F^*$ is opposite with some chamber containing $F^*$ or $F$ respectively. The opposition relation induces a (trivial or involutory) permutation $\tau^{\mathrm{op}}$ on the set $I$ of types of $\Delta$ such that if $F$ and $F^*$ are opposite flags of $\Delta$ then $\tau(F^*) = \tau^{\mathrm{op}}(\tau(F))$ (Tits \cite{Tits}). For the sake of completeness, we  recall that if $\Delta$ admits at least one irreducible component of type $A_n$ with $n > 1$, $D_n$ with $n$ odd, $E_6$ or $I_2(m)$ with $m$ odd then $\tau^{\mathrm{op}}$ is an involution, otherwise $\tau^{\mathrm{op}} = \mathrm{id}_I$ (Tits \cite{Tits}). 

Given a nonempty flag $A^*$ of $\Delta$, we say that a flag $F$ of $\Delta$ (possibly an element) is {\em far} from $A^*$ if $F$ is incident with at least one flag opposite to $A^*$. Let $\mathrm{Far}_\Delta(A^*)$ be the set of elements of $\Delta$ far from $A^*$, equipped with the type-function inherited from $\Delta$ and the incidence relation defined as follows: two element $x, y\in \mathrm{Far}_\Delta(A^*)$ are declared to be incident precisely when they are incident in $\Delta$ and the flag $\{x,y\}$ is far from $A^*$. Thus $\mathrm{Far}_\Delta(A^*)$ is a (possibly non-residually connected or even non-connected) geometry of rank $n$. The flags of $\mathrm{Far}_\Delta(A^*)$ are precisely the flags of $\Delta$ far from $A^*$. In particular, the flags of $\mathrm{Far}_\Delta(A^*)$ of type $\tau^{\mathrm{op}}(\tau(A^*))$ are the flags of $\Delta$ opposite to $A^*$.  

The following definition is needed in view of the next proposition. For two distinct types $i, j\in I$, let $g_{i,j}$, $s_i$ and $s_j$ be the gonality and the orders at $i$ and $j$ respectively of a $\{i,j\}$-residue of $\Delta$. (We recall that, since $\Delta$ is thick, the orders $s_i$ and $s_j$ do not depend on the choice of a particular $\{i,j\}$-residue.) 

\begin{prop}[{\rm Blok and Brouwer \cite{BB}}] \label{connessione}
Assume the following: 
\begin{itemize}
\item[$(*)$] For any two distinct types $i, j\in I$, if $i \in J$ then the triple $(g_{i,j},s_i,s_j)$ is different from either of $(6,2,2)$ and $(8,2,4)$. If $\{i,j\}\subseteq J$ then $(g_{i,j},s_i, s_j)$ is also different from either $(4,2,2)$ or $(6,3,3)$.
\end{itemize} 
Then the geometry $\mathrm{Far}_\Delta(A^*)$ is residually connected. 
\end{prop}

The next proposition is implicit in Ronan \cite[Chapter 6]{RonBook} (but see also Blok \cite[Lemmas 3.5 and 3.6]{Blok}). 

\begin{prop}\label{far residui}
Let $F$ be a flag of $\mathrm{Far}_\Delta(A^*)$. The residue $\mathrm{Res}_{\mathrm{Far}_\Delta(A^*)}(F)$ of $F$ in $\mathrm{Far}_\Delta(A^*)$ is the subgeometry of $\mathrm{Res}_\Delta(F)$ far from the projection $\mathrm{pr}_F(A^*)$ of $A$ onto $\mathrm{Res}_\Delta(F)$. In particular, $\mathrm{Res}_{\mathrm{Far}_\Delta(A^*)}(F) = \mathrm{Res}_\Delta(F)$ if and only if $\mathrm{pr}_F(A^*) = \emptyset$. We have $\mathrm{pr}_F(A^*) = \emptyset$ if and only if $\tau^{\mathrm{op}}(\tau(A))\subseteq \tau(F)$.
\end{prop}  

Let $F$ be a flag of $\mathrm{Far}_\Delta(A^*)$ of corank $2$. It readily follow from Proposition \ref{far residui} that $\mathrm{Res}_{\mathrm{Far}_\Delta(A^*)}(F)$ is a generalized digon if and only if $\mathrm{Res}_\Delta(F)$ is a generalized digon. Hence $\mathrm{Far}_\Delta(A^*)$ and $\Delta$ have the same diagram graph. 

\subsubsection{Shadows in $\mathrm{Far}_\Delta(A^*)$}\label{shadows in far}

Assume that $\mathrm{Far}_\Delta(A^*)$ is residually connected (compare Proposition \ref{connessione}). Let $\mathrm{op}(A^*)$ be the set of flags of $\Delta$ opposite to $A^*$. Let ${\cal F}(\mathrm{Far}_\Delta(A^*))$ be the set of nonempty flags of $\mathrm{Far}_\Delta(A^*)$. Put $J := \tau^{\mathrm{op}}(\tau(A^*))$ and assume that $J$ meets all components of the diagram of $\Delta$ non-trivially. For $F\in{\cal F}(\mathrm{Far}_\Delta(A^*))$, let $\mathrm{sh}_{J,A^*}(F) := \mathrm{sh}_J(F)\cap\mathrm{op}(A^*)$. Then the set 
\[\mathrm{Far}_\Delta(A^*)_J := \{\mathrm{sh}_{J,A^*}(F)\}_{F\in{\cal F}(\mathrm{Far}_\Delta(A^*))}\]
equipped with inclusion as the incidence relation, is a firm and residually connected geometry of rank $n$ (compare \cite[Chapter 5]{PasDG}). Properties quite similar to those that hold for $\Delta_J$ also hold for $\mathrm{Far}_\Delta(A^*)_J$. In particular, given a $J$-shadow $X$ in $\mathrm{Far}_\Delta(A^*)_J$, there exists a unique minimal flag $F_X\in{\cal F}(\mathrm{Far}_\Delta(A^*))$ such that $\mathrm{sh}_{J,A^*}(F_X) = X$. The $J$-shadow $\mathrm{sh}_J(F)$ of $F_X$ in $\Delta$ is the minimal element of $\Delta_J$ that contains $X$ and we have $X = \mathrm{sh}_J(F)\cap \mathrm{op}(A^*)$. For two $J$-shadows $X, Y \in \mathrm{Far}_\Delta(A^*)_J$ we have $X\subseteq Y$ if and only if $\mathrm{sh}_J(F_X)\subseteq \mathrm{sh}_J(F_Y)$ in $\Delta$. It follows that the function mapping a $J$-shadow $X\in\mathrm{Far}_\Delta(A)_J$ onto $\mathrm{sh}_J(F_X)$ yields an isomorphism from $\mathrm{Far}_\Delta(A^*)_J$ to the subgeometry of $\Delta_J$ induced on the set of $J$-shadows of $\Delta$ that meet $\mathrm{op}(A^*)$ non-trivially. 

Given a $J$-shadow $X \in \mathrm{Far}_\Delta(A^*)_J$ of $\mathrm{Far}_\Delta(A^*)$ of type greater than $1$ (namely $X$ is not a single $J$-flag) let  $F_X$ be defined as above. Then $J\not\subseteq\tau(F_X)$. Hence $\mathrm{pr}_{F_X}(A^*)\neq \emptyset$ and by the last claim of Proposition \ref{far residui}. The first part of Proposition \ref{far residui} implies the following, where $\mathrm{Res}_\Delta^J(F_X)$ is as in Subsection \ref{shadows residues sec}. 

\begin{co}\label{residui-shadow-far}
The elements of $\mathrm{Far}_\Delta(A^*)_J$ contained in $X$ form a geometry isomorphic to the $J\setminus\tau(F_X)$-shadow geometry of the subgeometry of $\mathrm{Res}^J_\Delta(F_X)$ far from $\mathrm{pr}_{F_X}(A^*)$.
\end{co} 

Let $\mathrm{Tr}^{\{1,2,3\}}(\mathrm{Far}_\Delta(A^*)_J)$ be the $\{1,2,3\}$-truncation of $\mathrm{Far}_\Delta(A^*)_J$, namely the subgeometry induced by $\mathrm{Far}_\Delta(A^*)_J$ on the set of elements of type $1$, $2$ or $3$. The next proposition is a special case of a well known property of truncations (see e.g. \cite[Theorem 1]{PasComo}).

\begin{prop}\label{truncation}
Suppose that for every element $X$ of $\mathrm{Far}_\Delta(A^*)_J$  of type $1$, $2$ or $3$ the residue of $X$ in $\mathrm{Far}_\Delta(A^*)_J$ is $2$-simply connected. Then the universal cover of $\mathrm{Tr}^{\{1,2,3\}}(\mathrm{Far}_\Delta(A^*)_J)$ is the $\{1,2,3\}$-truncation of the universal $2$-cover of $\mathrm{Far}_\Delta(A^*)_J$. In particular, $\mathrm{Tr}^{\{1,2,3\}}(\mathrm{Far}_\Delta(A^*)_J)$ is simply connected if and only if $\mathrm{Far}_\Delta(A^*)_J$ is $2$-simply connected.
\end{prop}

\section{Unipotent representations}\label{epsilonA-Sec}

Throughout this section $\Delta$ is a thick building of spherical type and rank $n\geq 3$ and $I$ is its set of types. To simplify our exposition we assume that $\Delta$ is irreducible, but what we are going to say in this section holds in the reducible case as well, provided that we assume that all irreducible components of $\Delta$ of rank $2$ are Moufang and all components of rank $1$ are projective lines. 

Let $G = \mathrm{Aut}(\Delta)$ be the group of all (type-preserving) automorphisms of $\Delta$ and $G^\infty$ the perfect core of $G$, namely the largest normal perfect subgroup of $G$. As $\Delta$ it is thick, irreducible and $n\geq 3$ by assumption,  the building $\Delta$ is Moufang (Tits \cite{Tits}).  Hence $G^\infty$ is generated by the root subgroups of $\Delta$ and it acts transitively on the set of pairs $(C,\Sigma)$ where $\Sigma$ is an apartment of $\Delta$ and $C$ a chamber of $\Sigma$.

For a nonempty flag $A^*$ of $\Delta$ we denote by $U_{A^*}$ the unipotent radical of the stabilizer $G_{A^*}$ of $A^*$ in $G$ and by $K_{A^*}$ the elementwise-stabilizer of $\mathrm{Res}_\Delta(A^*)$. We recall that $U_{A^*}$ is the subgroup of $G^\infty_{A^*} = G^\infty\cap G_{A^*}$ generated by the root subgroups $U_\alpha$ for $\alpha$ a root of (an apartment of) $\Delta$ such that $A^*$ is contained in $\alpha$ but not in the wall $\partial \alpha$ of $\alpha$ (see Ronan \cite[Chapter 6]{RonBook}). The group $U_{A^*}$ is the largest nilpotent subgroup of $K_{A^*}$.

\subsection{Transitivity properties of unipotent radicals}

The following is well known (see e.g. Ronan \cite[Chapter 6, exercise 17]{RonBook}). 

\begin{lemma}\label{regularity}
The group $U_{A^*}$ acts sharply transitively on the set $\mathrm{op}(A^*)$ of flags opposite to $A^*$. In particular, the identity is the unique element of $U_{A^*}$ that fixes a flag opposite to $A^*$. 
\end{lemma}

Let $J := \tau^{\mathrm{op}}(\tau(A^*))$. Given another flag $F$ of $\Delta$, let $\mathrm{op}_F(A^*)$ be the set of $J$-flags opposite to $A^*$ and incident with $F$. The next lemma is a generalization of Lemma \ref{regularity}.  

\begin{lemma}\label{unipotent}
Let $\mathrm{op}_F(A^*)\neq\emptyset$. Then $U_{A^*}\cap G_F$ acts sharply transitively on the set $\mathrm{op}_F(A^*)$
\end{lemma}
\pr  We may assume that $J\not\subseteq \tau(F)$, otherwise there is nothing to prove. By assumption, $F$ is incident with at least one flag opposite to $A^*$. Hence there exists at least one flag $F^*$ incident with $A^*$ and opposite to $F$. Put $B^* := F^*\cup A^*$ and for $A\in \mathrm{op}_F(A^*)$ let $B_A := A\cup F$. Then the flags $B^*$ and $B_A$ are opposite, as one can see by considering an apartment containing them both. Pick one $J$-flag $A_0\in \mathrm{op}_F(A^*)$ and let $\Sigma_0$ be an apartment containing both $B^*$ and $B_{A_0}$. Let $\alpha^+$ be a root of $\Sigma_0$ containing $B^*$ but with $B^*\cap \partial\alpha^+ = F^*$. Hence $A^*\not\subseteq \partial\alpha^+$, as $\tau(A^*)\not\subseteq\tau(F^*)$ since $J\not\subseteq\tau(F)$. Then $B_{A_0}$ is contained in the root $\alpha^-$ of $\Sigma_0$ opposite to $\alpha^+$ and $B_{A_0}\cap\partial\alpha^+ = F$. Moreover, the projection $\mathrm{pr}_F(A^*)$ of $A^*$ onto $\mathrm{Res}_\Delta(F)$ is contained in $\alpha^+$. Every member $A\in \mathrm{op}_F(A^*)$ belongs to an apartment $\Sigma$ such that the chambers of $\Sigma\cap\Sigma_0$ are precisely those of $\alpha^+$. An element $u\in U_{A^*}$ exists that fixes $\alpha^+$ element-wise and maps $\Sigma_0$ onto $\Sigma$. So, $u$ maps $A_0$ onto $A$ and fixes $F$. The transitivity of $U_{A^*}$ on $\mathrm{op}_F(A^*)$ is proved. By the second claim of Lemma \ref{regularity}, the action of $U_{A^*}$ on $\mathrm{op}_F{A^*}$ is sharply transitive. \eop

\begin{co}\label{unipotent-bis}
With $J$, $A^*$ and $F$ as in Lemma {\rm \ref{unipotent}} and ${\cal D}^J_F$ as in Subsection {\rm \ref{shadows}}, let $B$ be a flag of cotype ${\cal D}_{F}^J$ containing $F$. Then $U_{A^*}\cap G_F = U_{A^*}\cap G_B$. 
\end{co}
\pr Clearly, $U_{A^*}\cap G_B\leq U_{A^*}\cap G_F$, since $G_B\leq G_F$. So, we only must prove that $U_{A^*}\cap G_B\geq U_{A^*}\cap G_F$.

We have $\mathrm{sh}_J(F) = \mathrm{sh}_J(B)$ (Tits \cite[Chapter 12]{Tits}; also \cite[Chapter 5]{PasDG}). For $u\in U_{A^*}\cap G_F$ let $A' := u(A)$. By Lemma \ref{unipotent} and since $A$ and $A'$ belong to 
$\mathrm{sh}_J(F) = \mathrm{sh}_J(B)$ and are opposite to $A^*$, there exists an element $v\in U_{A^*}\cap G_B$ such that $v(A) = A'$. So, $u$ and $v$ are elements of $U_{A^*}$ mapping $A$ onto the same $J$-flag $A'$. It follows that $v = u$, by the second claim of Lemma \ref{regularity}. Hence $u \in U_{A^*}\cap G_B$. Therefore $U_{A^*}\cap G_F\leq U_{A^*}\cap G_B$ . \eop

\subsection{The embedding $\varepsilon_\Delta^J$}\label{sezione-epsilonA}

For $\emptyset \neq J\subseteq I$ and a $J$-flag $A$, let $A^*$ be a flag of $\Delta$ opposite to $A$. Let $\Delta_J(A)$ the local geometry of $\Delta_J$ at $A$ (Subsection \ref{local}). We define an embedding $\varepsilon_\Delta^A:\Delta_J(A)\rightarrow U_{A^*}$ of $\Delta_J(A)$ into the unipotent radical $U_{A^*}$ of $G_{A^*}$ as follows: for every element (point or line) $X$ of $\Delta_J(A)$ we set 
\[\varepsilon^A_\Delta(X) := U_{A^*}\cap G_{F_X} ~ (= U_{A^*}\cap G^\infty_{F_X}),\]
where $F_X$ is the minimal element of ${\cal F}_X$, namely the minimal flag $F$ of $\Delta$ such that $\mathrm{sh}_J(F) = X$ (Subsection \ref{shadows}). We call $\varepsilon^A_\Delta$ the {\em unipotent representation} of $\Delta_J(A)$. 

The next theorem states that $\varepsilon^A_\Delta$ is indeed an embedding, except possibly in the exceptional cases excluded by $(*)$ of Proposition \ref{connessione}.  

\begin{theo}\label{epsilonA}
Both the following hold. 

\begin{itemize}
\item[$(1)$] If $X_1$ and $X_2$ are distinct points of $\Delta_J(A)$, then $\varepsilon^A_\Delta(X_1)\cap \varepsilon^A_\Delta(X_2) = 1$. 
\item[$(2)$] For a point $X\in{\cal P}_A$ and a line $Y\in{\cal L}_A$, we have $\varepsilon^A_\Delta(X)\leq \varepsilon^A_\Delta(Y)$ if and only if $X \in P(Y)$.
\end{itemize}
If moreover $(*)$ of Proposition {\rm \ref{connessione}} holds then the following also hold:

\begin{itemize}
\item[$(3)$] $\varepsilon^A_\Delta(Y) = \langle \varepsilon^A_\Delta(X)\rangle_{X\in P(Y)}$ for any line $Y\in{\cal L}_A$.
\item[$(4)$] $U_{A^*} =\langle \varepsilon^A_\Delta(X)\rangle_{X\in{\cal P}_A}$.
\end{itemize}
\end{theo}
\pr 
Let $X_1$ and $X_2$ be distinct points of $\Delta_J(A)$. Then $X_1\cap X_2 = \{A\}$. Hence $U_{A^*}\cap G_{F_{X_1}}\cap G_{F_{X_2}}$ stabilizes $A$. Therefore $U_{A^*}\cap G_{F_{X_1}}\cap G_{F_{X_2}} = 1$, since the identity is the only element of $U_{A^*}$ that can fix a $J$-flag opposite $A^*$. Claim $(1)$ is proved. 

We shall now prove claim (2). Let $X$ and $Y$ be a point and a line of $\Delta_J(A)$. Assume firstly that $X\subseteq Y$.  Then ${\cal D}_X\subseteq {\cal D}_Y$ and $F_X\cup F_Y$ is a flag. Let $C$ be a flag of $\Delta$ of cotype ${\cal D}_X$ containing $F_X\cup F_Y$. Then $U_{A^*}\cap G_C\leq U_{A^*}\cap F_Y$. Moreover $U_{A^*}\cap G_C = U_{A^*}\cap G_{F_X}$ by Corollary \ref{unipotent-bis}. Therefore $U_{A^*}\cap G_{F_X}\leq U_{A^*}\cap G_{F_Y}$. 

Conversely, let $U_{A^*}\cap G_{F_X}\leq U_{A^*}\cap G_{F_Y}$. Let $\overline{X} := \mathrm{op}_{F_X}(A^*) \subseteq X$ and $\overline{Y} := \mathrm{op}_{F_Y}(A^*)\subseteq Y$. By Lemma \ref{unipotent} the sets $\overline{X}$ and $\overline{Y}$ are orbits of $U_{A^*}\cap G_{F_X}$ and $U_{A^*}\cap G_{F_Y}$ respectively. It follows that $\overline{X}\subseteq \overline{Y}$, since $A\in \overline{X}\cap \overline{Y}$ and $U_{A^*}\cap G_{F_X}\leq U_{A^*}\cap G_{F_Y}$ by assumption. Hence $\overline{X}\subseteq Y$. On the other hand, as $|X| > 1$ and $\Delta$ is thick, we have $|\overline{X}| > 1$. Therefore $|X\cap Y| > 1$. This forces $X\in P(Y)$. 

Turning to $(3)$, suppose that $(*)$ of Proposition \ref{connessione} holds. Let $Y\in{\cal L}_A$. By Proposition \ref{shadows-residues}, the $J\setminus{\cal D}_Y$-shadow geometry $\mathrm{Res}_\Delta^J(F_Y)_{J\setminus{\cal D}_Y}$ of $\mathrm{Res}_\Delta^J(F_Y)$ is isomorphic to $\mathrm{Res}_{\Delta_J}^-(Y)$. The function that maps every flag $F$ of $\mathrm{Res}_\Delta^J(F_Y)$ onto the flag $F\cup (F_Y\cap A)$ is indeed an isomorphism from $\mathrm{Res}_{\Delta_J}^-(Y)$ to $\mathrm{Res}_\Delta^J(F_Y)_{J\setminus{\cal D}_Y}$. The members of $\mathrm{op}_{F_Y}(A^*)$ and the elements of $\mathrm{Res}_{\Delta_J}^-(Y)$ incident with at least two members of $\mathrm{op}_{F_Y}(A^*)$ form a subgeometry $\mathrm{Res}^-_{\Delta_J, A^*}(Y)$ of $\mathrm{Res}_{\Delta_J}^-(Y)$. We take $\mathrm{op}_{F_Y}(A^*)$ as the point-set of $\mathrm{Res}^-_{\Delta_J, A^*}(Y)$, the remaining elements of  $\mathrm{Res}^-_{\Delta_J, A^*}(Y)$ being the lines of  $\mathrm{Res}^-_{\Delta_J, A^*}(Y)$. Thus, $A$ is a point of $\mathrm{Res}^-_{\Delta_J, A^*}(Y)$.

By Corollary \ref{residui-shadow-far}, the geometry $\mathrm{Res}^-_{\Delta_J, A^*}(Y)$ is isomorphic to the subgeometry of $\mathrm{Res}^J_{\Delta}(F_Y)$ far from $B$. By Proposition \ref{connessione} and since $\Delta$ and $\tau^{\mathrm{op}}(J)$ are not as in any of the exceptional examples mentioned in that proposition, the subgeometry of $\mathrm{Res}^J_{\Delta}(F_Y)$ far from $B$ is connected. Hence $\mathrm{Res}^-_{\Delta_J, A^*}(Y)$ is connected as well. Let $u\in U_{A^*}\cap G_{F_Y}$ and $A' := u(A)$. Clearly, $A'$ is a point of $\mathrm{Res}^-_{\Delta_J, A^*}(Y)$. Let $d$ be the distance from $A$ to $A'$ in the collinearity graph of $\mathrm{Res}^-_{\Delta_J, A^*}(Y)$. By induction on $d$, we shall prove that $u\in \langle U_{A^*}\cap G_{F_X}\rangle_{X\in P(Y)}$. 

If $d = 0$ there is nothing to prove. Indeed in this case $A' = A$, hence $u = 1$. Let $d = 1$ and let $X$ be the line of  $\mathrm{Res}^-_{\Delta_J, A^*}(Y)$ incident with both $A$ and $A'$.  By Lemma \ref{unipotent}, the subgroup $U_{A^*}\cap G_{F_X}$ contains an element $v$ such that $v(A) = A' = u(A)$. It follows that $u = v$, hence $u \in U_{A^*}\cap G_{F_X}$. 

Let $d > 1$ and let $A = A_0, A_1,..., A_{d-1}, A_d = A'$ be a minimal path from $A$ to $A'$ in the collinearity graph of $\mathrm{Res}^-_{\Delta_J,A^*}(Y)$. Let $u_1, u_2\in U_{A^*}\cap G_{F_Y}$ map $A$ onto $A_{d-1}$ and $A_{d-1}$ onto $A'$, respectively. Elements $u_1$ and $u_2$ with these properties exist by Lemma \ref{unipotent}. Moreover $u = u_2u_1$, since $u(A) = A' = u_2(u_1(A))$. By the inductive hypothesis, $u_1\in  \langle U_{A^*}\cap G_{F_X}\rangle_{X\in P(Y)}$. Moreover, $A'' := u_1^{-1}(A')$ is collinear with $A$ in $\mathrm{Res}^-_{\Delta_J, A^*}(Y)$ and $u_1^{-1}u_2u_1$ maps $A$ onto $A'$. Hence $u_1^{-1}u_2u_1\in U_{A^*}\cap G_{F_X}$ for some $X\in P(Y)$, by the previous paragraph. It follows that $u = u_1\cdot u_1^{-1}u_2u_1$ belongs to
$\langle U_{A^*}\cap G_{F_X}\rangle_{X\in P(Y)}$. Claim $(3)$ is proved.

Claim $(4)$ can be proved by an argument very similar to the one used for $(3)$. We leave the details for the reader. As in the proof of $(3)$, hypothesis $(*)$ is needed to ensure that the point-line geometry formed by the $J$-flags opposite to $A^*$ and the lines of $\Delta_J$ incident with at least two such $J$-flags is connected. \eop  

\bigskip

Henceforth we always assume that $(*)$ of Proposition \ref{connessione} holds. 

\begin{theo}\label{1.5}
The expansion $\mathrm{Ex}_{\varepsilon^A_\Delta}(\Delta_J(A))$ of $\Delta_J(A)$ in $U_{A^*}$ by $\varepsilon^A_\Delta$ is isomorphic to the $\{1,2,3\}$-truncation $\mathrm{Tr}^{\{1,2,3\}}(\mathrm{Far}_\Delta(A^*)_J)$ of the $J$-shadow geometry $\mathrm{Far}_\Delta(A^*)_J$ of $\mathrm{Far}_\Delta(A^*)$. 
\end{theo}  
\pr 
As $U_{A^*}$ is sharply transitive on $\mathrm{op}(A^*)$, a bijection can be established between $\mathrm{op}(A^*)$ and $U_{A^*}$ which maps every $J$-flag $A'\in \mathrm{op}(A^*)$ onto the unique element $u\in U_{A^*}$ such that $u(A) = A'$. This bijection naturally extends to an isomorphism from $\mathrm{Tr}^{\{1,2,3\}}(\mathrm{Far}_\Delta(A^*)_J)$ to $\mathrm{Ex}({\varepsilon^A_\Delta})$ as follows. Let $\overline{X}$ be an element of $\mathrm{Far}_\Delta(A^*)_J$ of type $i\in\{2,3\}$. Then $\overline{X} = X\cap \mathrm{op}(A^*)$ for a unique element $X$ of $\Delta_J$ of type $\tau_J(X) = i$ (see Subsection \ref{shadows in far}). Let $u\in U_{A^*}$ map $A$ onto a $J$-flag incident to $F_X$. Then $\overline{X}$ is mapped onto the coset $u\cdot U_{A^*}\cap G_{u^{-1}(F_X)}$. It is straighforward to check that the mapping defined in this way is indeed an isomorphism.  \eop

\bigskip

 By combining Theorem \ref{1.5} with Proposition \ref{expansion} we immediately obtain the following:

\begin{co}\label{espansione-epsilonA-co}
The embedding $\varepsilon_\Delta^A$ is dominant if and only if the geometry $\mathrm{Tr}^{\{1,2,3\}}(\mathrm{Far}_\Delta(A^*)_J)$ is simply connected. 
\end{co}  

By Lemma \ref{regularity}, the stabilizer $G_A$ of $A$ in $G$ is a semi-direct product of $U_A$ and $G_{A,A^*} := G_A\cap G_{A^*}$. As $A^*$ and $A$ are opposite, every chamber containing $A^*$ is the gate of ${\cal C}(A^*)$ to a unique chamber containing $A$ (see Subsection \ref{proj}). Hence the intersection $K_A\cap G_{A,A^*}$ is contained in $K_{A^*}$. By symmetry, $K_{A}^*\cap G_{A,A^*}\leq K_{A}$. Therefore $K_A\cap G_{A,A^*} = K_{A^*}\cap G_{A,A^*}$. Assume the following:

\begin{itemize}
\item[$(\mathrm{A})$] Every automorphism $\alpha$ of $\Delta_J(A)$ is induced by an element $\alpha_A$ of $G_{A,A^*}$ (uniquely determined modulo $K_A\cap G_{A,A^*}$). 
\end{itemize}
Let $\alpha_A$ be as in $(\mathrm{A})$. Then $\alpha_A$ normalizes $U_{A^*}$, since $U_{A^*}$ is a characteristic subgroup of $K_{A^*}$ and the latter is normal in $G_{A^*}$, which contains $\alpha_A$. Let $\alpha^*$ be the automorphism of $U_{A^*}$ induced by $\alpha_A$. It is not so difficult to prove that $\varepsilon_\Delta^A( \alpha(X)) = \alpha^*(\varepsilon_\Delta^A(X))$ for every element $X$ of $\Delta_J(A)$. In other words, $\alpha^*$ is a lifting of $\alpha$. Therefore,

\begin{theo}\label{epsilonA-omogeneo}
Let $(\mathrm{A})$ hold. Then the embedding $\varepsilon_\Delta^A$ is homogeneous.
\end{theo}

Still assuming that $(\mathrm{A})$ holds, let $\Gamma$ be a point-line geometry isomorphic to $\Delta_J(A)$. Given an isomorphism $\gamma:\Gamma\rightarrow\Delta_J(A)$, put $\varepsilon_{\gamma,\Delta}^J := \varepsilon_\Delta^A\cdot\gamma$. Then $\varepsilon_{\gamma,\Delta}^J$ is an embedding of $\Gamma$ in $U_{A^*}$. By Theorem \ref{epsilonA-omogeneo}, if $\gamma_1$ and $\gamma_2$ are two isomorphisms from $\Gamma$ to $\Delta_J(A)$, then $\varepsilon_{\gamma_1, \Delta}^J\cong \varepsilon_{\gamma_2,\Delta}^J$. In view of this, we feel free to write $\varepsilon_\Delta^J$ for $\varepsilon_{\gamma,\Delta}^J$, dropping the subscript $\gamma$. We call $\varepsilon_\Delta^J$ the {\em unipotent representation} of $\Gamma$ {\em of type} $(\Delta, J)$. 

When $J$ is a singleton, say $J = \{i\}$, we write $\varepsilon^i_\Delta$ instead of $\varepsilon^{\{i\}}_\Delta$. 

\bigskip

\noindent
{\bf Remark.} It is not difficult to find examples where $(\mathrm{A})$ fails to hold. For instance, this is the case when $J = I$, but it is likely that less trivial counterexamples exist. 

\subsection{Representations induced on residues}\label{indotto sezione} 

Assume that $(*)$ of Proposition \ref{connessione} holds. Let $F$ be a flag incident with $A$ but not containing $A$ and such that $\mathrm{Res}_\Delta^J(F)$ (see Subsection \ref{shadows residues sec}) has rank at least $3$. Put $\Delta_F := \mathrm{Res}_\Delta^J(F)$, $J_F := J\setminus \tau(F)$, $A_F = A\setminus F$ and $S_F := \mathrm{sh}_J(F)$. The $J_F$-shadow geometry $\Delta_{F,J_F}$ of $\Delta_F$ is isomorphic $\mathrm{Res}^-_{\Delta_J}(S_F)$ (Proposition \ref{shadows-residues}). Let ${\cal X}_{A,F}$ be the set of elements of $\Delta_J(A)$ contained in $S_F$ and let $\Delta_J(A)_F$ be the subgeometry induced by $\Delta_J(A)$ on ${\cal X}_{A,F}$. The isomorphism $\Delta_{F,J_F}\cong \mathrm{Res}^-_{\Delta_J}(S_F)$ induces an isomorphism between the local geometry $\Delta_{F,J_F}(A_F)$ of $\Delta_{F,J_F}$ at $A_F$ and the geometry of $\Delta_{J}(A)_F$, mapping every element $X\in{\cal X}_{A,F}$ onto the $J_F$-shadow $\overline{X} := \{Y\setminus F\}_{Y\in X}$. 

Clearly, $(*)$ is inherited by the pair $(\Delta_F, J_F)$. Hence the unipotent representation $\varepsilon^{A_F}_{\Delta_F}$ of $\Delta_{F, J_F}(A_F)$ is an embedding. 

\begin{lemma}\label{non fa differenza}
We have $U_{A^*}\cap G_{F_X} = U_{A^*}\cap G_{F_X}\cap G_F$ for any $X\in{\cal X}_{A,F}$.
\end{lemma}
\pr We have $\mathrm{sh}_J(F_X) = X \subseteq S_F = \mathrm{sh}_J(F)$. Hence $\mathrm{sh}_J(F_X\cup F) = X$. By this equality and Lemma \ref{unipotent} each of $U_{A^*}\cap G_{F_X}$ and $U_{A^*}\cap G_{F_X\cup F}$ acts sharply transitively on $X$. However $U_{A^*}\cap G_{F_X\cup F}\leq U_{A^*}\cap G_{F_X}$. It follows that $U_{A^*}\cap G_{F_X \cup F} = U_{A^*}\cap G_{F_X}$. \eop

\bigskip 

Let $\varepsilon_{\Delta,F}^A$ be the restriction of $\varepsilon_\Delta^A$ to $\Delta_J(A)_F$. Then $\varepsilon_{\Delta,F}^A$ is an embedding of $\Delta_J(A)_F$ in the subgroup $U^F_{A^*}$ of $U_{A^*}$ generated by the subgroups $U_{A^*}\cap G_{F_X}$ for $X\in{\cal X}_{A,F}$. Note that $U^F_{A^*}\leq U_{A^*,F}$, by Lemma \ref{non fa differenza}. In view of the isomorphism $\Delta_{F,J_F}(A_F)\cong \Delta_J(A)_F$ we can also regard $\varepsilon_{\Delta,F}^A$ as an embedding of $\Delta_{F,J_F}(A_F)$ in the subgroup $U^F_{A^*}$ of $U_{A^*,F}$.  

\begin{theo}\label{unipotente indotto}
We have $U^F_{A^*} = U_{A^*,F}$ and $\varepsilon_{\Delta,F}^A\cong \varepsilon_{\Delta_F}^{A_F}$.  
\end{theo}
\pr As remarked in Subsection \ref{shadows residues sec}, we may assume that $\tau(F) = {\cal D}^J_F$. So, $\Delta_F = \mathrm{Res}_\Delta(F)$. We have $U_{A^*}\cap K_F = 1$ by the second claim of Lemma \ref{unipotent}. Hence $U_{A^*,F}$ can be regarded as a subgroup of $\overline{G}_F := (G_FK_F)/K_F$. Clearly, $U_{A^*,F}$ stabilizes $B := \mathrm{pr}_F(A^*)$. Let $\overline{U}_B$ be the unipotent radical of the stabilizer $\overline{G}_{F,B}$ of $B$ in $\overline{G}_F$. 

Let $g\in G^\infty_{F\cup B}$. Then $\mathrm{pr}_F(A^*) = \mathrm{pr}_F(g(A^*)) = B$. Let $A$ be a $J$-flag incident to $F$ and opposite to $A^*$. By Corollary \ref{residui-shadow-far}, the flags $A\setminus F$ and $B$ are opposite in $\mathrm{Res}_\Delta(F)$. Let $U_{A\cup F}$ and $\overline{U}_{A\setminus F}$ be the unipotent radicals of $G^\infty_{A\cup F}$ and $\overline{G}_{F, A\setminus F}$ respectively. Then 
\begin{equation}\label{unipotent eq}
\frac{U_{A\cup F}K_F}{K_F} = \overline{U}_{A\setminus F}.
\end{equation}
Given a flag $F^*$ of type $\tau^{\mathrm{op}}(\tau(F))$ incident with $A^*$ and such that $F^*\cup A^*$ is opposite to $F\cup A$, the flag $g(F^*\cup A^*)$ is opposite to $F\cup A$. (This follows from Corollary \ref{residui-shadow-far}, noticing that $\mathrm{pr}_F(F^*\cup A^*) = \mathrm{pr}_F(A^*) = B$ and recalling that $B$ and $A\setminus F$ are opposite.) Hence, by Lemma \ref{unipotent}, there exists an element $u\in U_{A\cup F}$ mapping $g(F^*\cup A^*)$ back to $F^*\cup A^*$. However $u$ stabilizes both $A\setminus F$ and $B = \mathrm{pr}_F(F^*\cup A^*) = \mathrm{pr}_F(g(F^*\cup A^*)$. By equality (\ref{unipotent eq}), we obtain that $u \in K_F$. Therefore $g\in G^\infty_{F,A^*}K_F$, where $G^\infty_{F,A^*} := G^\infty_F\cap G^\infty_{A^*}$. As a consequence, $(G^\infty_{F,A^*}K_F)/K_F = \overline{G}_{F, B}$. It follows that $(U_{A^*, F}K_F)/K_F$ is normalized by $\overline{G}_{F,B}$. Thus, $(U_{A^*,F}K_F)/K_F$ is a normal nilpotent subgroup of $\overline{G}_{F,B}$. It is contained in the largest normal nilpotent subgroup of $\overline{G}_{F,B}$, namely $\overline{U}_B$. Therefore $U_{A^*,F}$  ($\cong (U_{A^*,F}K_F)/K_F$ since $U_{A^*}\cap K_F = 1$) can be regarded as a subgroup of $\overline{U}_B$. By Lemma \ref{unipotent}, the group $U_{A^*,F}$ acts sharply transitively on the set $\mathrm{op}_F(A^*)$ of $J$-flags incident to $F$ and opposite to $A^*$. Consequently, $U_{A^*,F} = \overline{U}_B$.

The isomorphism $\varepsilon_{\Delta,F}^A\cong \varepsilon_{\Delta_F}^{A_F}$ as well as the equality $U^F_{A^*} = U_{A^*,F}$ immediately follow from the equality $U_{A^*,F} = \overline{U}_B$ and Lemma \ref{non fa differenza}. \eop

\subsection{The abelianization $\varepsilon^J_{\Delta,\mathrm{ab}}$ of $\varepsilon_\Delta^J$}\label{epsilonA-quotient-section}

We still assume that $(*)$ of Proposition \ref{connessione} holds. Let $U'_{A^*}$ be the derived subgroup of $U_{A^*}$. Clearly $\varepsilon^A_\Delta$ admits an abelian quotient if and only if $U'_{A^*}$ defines a quotient of $\varepsilon^A_\Delta$, namely it satisfies the following conditions (compare $(\mathrm{Q}1)$ and $(\mathrm{Q}2)$ of Subsection \ref{morphisms quotients}): 

\begin{itemize}
\item[$(\mathrm{B}1)$] $U'_{A^*}\cap G_{F_X} = 1$ for every element $X$ of $\Delta_J(A)$. 
\item[$(\mathrm{B}2)$] $((U_{A^*}\cap G_{F_X})U'_{A^*})\cap((U_{A^*}\cap G_{F_{Y}})U'_{A^*}) = U'_{A^*}$ for any two non-collinear points $X, Y \in{\cal P}_A$. 
\end{itemize}

Assuming that $U'_{A^*}$ defines a quotient of $\varepsilon^A_\Delta$ we denote the quotient $\varepsilon^A_\Delta/U'_{A^*}$ by the symbol $\varepsilon^A_{\Delta,\mathrm{ab}}$ and we call it the {\em abelianization} of $\varepsilon^A_\Delta$. 

Clearly, $\varepsilon^A_{\Delta,\mathrm{ab}}$ is the largest abelian quotient of $\varepsilon^A_\Delta$. Consequently, if $\varepsilon^A_\Delta$ is dominant then $\varepsilon^A_{\Delta,\mathrm{ab}}$ is its own abelian hull. Suppose moreover that the following holds. 

\begin{itemize}
\item[$(\mathrm{B}3)$] The factor group $U_{A^*}/U'_{A^*}$ is the additive group of a vector space $V$ and $((U_{A^*}\cap G_{F_X})U'_{A^*})/U'_{A^*}$ is (the additive group of) a vector subspace of $V$, for every element $X$ of $\Delta_J(A)$. 
\end{itemize}
Then $\varepsilon_{\Delta,\mathrm{ab}}^A$ is linear. It is indeed the largest linear quotient of $\varepsilon_\Delta^A$. So, if the latter is dominant then $\varepsilon^A_{\Delta,\mathrm{ab}}$ is linearly dominant. If the following also holds besides $(\mathrm{B}3)$, then $\varepsilon^A_{\Delta,\mathrm{ab}}$ is projective.

\begin{itemize}
\item[$(\mathrm{B}4)$] The vector space $((U_{A^*}\cap G_{F_X})U'_{A^*})/U'_{A^*}$ has dimension $1$ or $2$ according to whether $X$ is either a point or a line of $\Delta_J(A)$. 
\end{itemize}
Finally, let $\Gamma\cong \Delta_J(A)$ and suppose that $(\mathrm{A})$ of Subsection \ref{sezione-epsilonA} holds. We put $\varepsilon^J_{\Delta,\mathrm{ab}} = \varepsilon_{\Delta,\mathrm{ab}}^A\cdot \gamma$ for a given isomorphism from $\Gamma$ to $\Delta_J(A)$, no matter which, and we call $\varepsilon^J_{\Delta,\mathrm{ab}}$ the {\em abelianization} of $\varepsilon^J_\Delta$. 

\subsection{A few examples related to $\mathrm{PG}(n-1,\mathbb{K})$}\label{examples} 

In order to give the reader a more concrete idea of what unipotent representations look like, we shall now discuss a few examples. In each of them $\Gamma = A_{n-1,1}(\mathbb{K})$, the $1$-shadow geometry of the building $A_{n-1}(\mathbb{K})$ (notation as in Section \ref{notation}). Namely, $\Gamma$ is the projective space $\mathrm{PG}(n-1,\mathbb{K})$. We assume $n\geq 3$. 

Let $\Gamma$ be as above. Then $\Gamma$ admits an obvious unipotent representation $\varepsilon_\Delta^1$ in the building $\Delta = A_n(\mathbb{K})$, isomorphic to the natural embedding of $\Gamma$ in $V(n,\mathbb{K})$.

Assume that $\mathbb{K}$ is isomorphic to its opposite $\mathbb{K}^{\mathrm{op}}$. Then $\Gamma$ admits unipotent representations different from $\varepsilon_{A_n(\mathbb{K})}^1$. Indeed in this case there exists at least one thick building $\Delta$ of type $C_n$ such that the maximal singular subspaces of the polar space $\Delta_1$ associated to $\Delta$ are isomorphic to $\mathrm{PG}(n-1,\mathbb{K}^{\mathrm{op}})$. Let $A$ and $A^*$ be two opposite $n$-elements of $\Delta$. Then $\Gamma\cong \Delta_n(A)$ and $\varepsilon^n_\Delta$ is the unipotent representation of $\Gamma$ in the unipotent radical $U_{A^*}$ of the stabilizer of $A^*$ in $\mathrm{Aut}(\Delta)$. Hypothesis $(\mathrm{A})$ of Subsection \ref{sezione-epsilonA} holds. Thus $\varepsilon_\Delta^n$ is well defined, namely it does not depend on the choice of a particular isomorphism from $\Gamma$ to $\Delta_n(A)$. Property $(*)$ of Proposition \ref{connessione} also holds. Hence $\varepsilon^n_\Delta$ is indeed an embedding.

We shall examine five of these representations. In each of them $\mathbb{K}$ is commutative. Explicitly, we take $\Delta$ equal to $C_n(\mathbb{K})$,  $B_n(\mathbb{K})$, ${^2}A_{2n-1}(\mathbb{F})$, ${^2}A_{2n}(\mathbb{F})$ and ${^2}D_{n+1}(\mathbb{K})$ respectively (notation as in Subsection \ref{notation-Dynkin}).  

In each of the cases that we are going to examine either $U_A^*$ is abelian or $U'_{A^*}$ does not define a quotient of $\varepsilon_\Delta^n$. This fact depends on the following quite general property. Let $\varepsilon:\Gamma\rightarrow U$ be an embedding of $\Gamma$ in a group $U$. Since any two points of $\Gamma$ are joined by a line, if $\varepsilon(L)$ is commutative for any line $L$ of $\Gamma$ then $U$ is commutative (compare \cite[Lemma 5.7]{PasEE}). In other words, either $U'\cap \varepsilon(L)\neq 1$ for some line $L$ (hence $U'$ does not define a quotient of $\varepsilon$) or $U' = 1$. 

\subsubsection{The case $\Delta = C_n(\mathbb{K})$ (quadratic veronesean embeddings)}\label{3.5.1}

Let  $\Delta = C_n(\mathbb{K})$. The group $U_{A^*}$ is isomorphic to the additive group of the vector space $\mathrm{S}_n(\mathbb{K})$ of $n\times n$ symmetric matrices over $\mathbb{K}$. Regarding $U_{a^*}$ as a copy of $\mathrm{S}_n(\mathbb{K})$, the embedding $\varepsilon_\Delta^n$ is $\mathbb{K}$-linear but non-projective, with $\mathrm{dim}(\varepsilon_\Delta^n(X)) = 1$ or $3$ according to whether $X$ is a point or a line of $\Gamma$. In fact $\varepsilon^n_\Delta$ is isomorphic to the quadratic veronesean embedding of $\Gamma$ in $V({{n+1}\choose 2}, \mathbb{K})\cong U_{A^*}$ (see \cite[Section 6.1]{PasEE}). So, $\{\varepsilon_\Delta^n(p)\}_{p\in L}$ is a conic in the projective plane $\mathrm{PG}(\varepsilon^n_\Delta(L))$, for every line $L$ of $\Gamma$. 

\subsubsection{The case $\Delta = B_n(\mathbb{K})$}\label{BnK}

Let $\Delta = B_n(\mathbb{K})$. Let $V = V(n,\mathbb{K})$. The elements of $U_{A^*}$ can be regarded as pairs $(S+v^Tv,v)$ with $v = (v_i)_{i=1}^n\in V$ (regarded as a row-vector) and $S$ a skew-symmetric matrix of order $n$. The elements of $U_{A^*}$ of the form $(S,0)$ form a subgroup $W$ stabilized by $G_{A^*}$ and $U_{A^*}/W$ is (the additive group of) a copy of $V$. Condition $(\mathrm{Q}1)$ of Subsection \ref{morphisms quotients} fails to hold for $W$. Hence $W$ does not define a quotient of $\varepsilon^n_\Delta$. Nevertheless, regarding $U_{A^*}/W$ as a copy of $V$, the function that maps every element $X\in\Gamma$ onto $(\varepsilon_\Delta^n(X)W)/W$ is a projective embedding of $\Gamma$ in $U_{A^*}/W\cong V$, isomorphic to the natural embedding of $\Gamma$ in $V$.   

If $\mathrm{char}(\mathbb{K})\neq 2$ then $U'_{A^*} = W$. In this case $U'_{A^*}$ does not define a quotient of $\varepsilon_\Delta^n$, namely $\varepsilon_\Delta^n$ does not admit any abelian quotient.  

Let $\mathrm{char}(\mathbb{K}) = 2$. Then the projection $\pi: (S+v^Tv,v)\mapsto S+v^Tv$ yields an injective morphism from $U_{A^*}$ to the additive group of $n\times n$ symmetric matrices over $\mathbb{K}$. Hence $U_{A^*}$ is abelian. Note that $\pi$ is surjective if and only if $\mathbb{K}$ is perfect. When $\mathbb{K}$ is perfect  then $B_n(\mathbb{K}) \cong C_n(\mathbb{K})$. In this case $\varepsilon_\Delta^n$ is the quadratic veronesean embedding of $\Gamma$ (see Subsection \ref{3.5.1}), hence it is $\mathbb{K}$-linear. Suppose that $\mathbb{K}$ is non-perfect. Then we can put a $\mathbb{K}$-vector space structure $V'$ on $U_{A^*}$ by setting $t\cdot (S+v^Tv,v) := (t^2(S+v^Tv), tv)$ for every $t\in \mathbb{K}$ and every pair $(S+vv^T,v)\in U_{A^*}$. It is not difficult to see that
\[\mathrm{dim}(V') = n + {{n-1}\choose 2}\cdot[\mathbb{K}:\mathbb{F}],\]
where $\mathbb{F}$ is the subfield of square elements of $\mathbb{K}$. (Note that $[\mathbb{K}:\mathbb{F}]$ might be infinite.) With $U_{A^*}$ regarded as a vector space in this way, the embedding $\varepsilon_\Delta^n$ is linear but non-projective. We have $\mathrm{dim}(\varepsilon_\Delta^n(p)) = 1$ for every point $p$ of $\Gamma$ while if $L$ is a line then $\mathrm{dim}(\varepsilon_\Delta^n(L)) = 2 + [\mathbb{K}:\mathbb{F}]$. 
 
\subsubsection{The case $\Delta = {^2}A_{2n-1}(\mathbb{F})$ (hermitian veronesean embeddings)}

In this case $\mathbb{K}$ is a quadratic extension of a subfield $\mathbb{F}$ and $\Delta = {^2}A_{2n-1}(\mathbb{F})$. Using an anti-hermitian form to define $\Delta$, one can see that $U_{A^*}$ is isomorphic to the additive group of the $\mathbb{F}$-vector space $\mathrm{H}_n(\mathbb{K})$ of $n\times n$ hermitian matrices over $\mathbb{K}$. Regarding $U_{a^*}$ as a copy of $\mathrm{H}_n(\mathbb{K})$, the embedding $\varepsilon_\Delta^n$ is $\mathbb{F}$-linear but non-projective, with $\mathrm{dim}(\varepsilon_\Delta^n(X)) = 1$ or $4$ according to whether $X$ is a point or a line of $\Gamma$. In fact $\varepsilon^n_\Delta$ is isomorphic to the hermitian veronesean embedding of $\Gamma$ in $V(n^2, \mathbb{F})\cong U_{A^*}$ (see \cite[Section 6.2]{PasEE}, where hermitian veronesean embeddings are called twisted tensor embeddings). Hence $\{\varepsilon_\Delta^n(p)\}_{p\in L}$ is an elliptic quadric in the $3$-dimensional $\mathbb{F}$-projective space $\mathrm{PG}(\varepsilon_\Delta^n(L))$, for every line $L$ of $\Gamma$. 

\subsubsection{The case $\Delta = {^2}A_{2n}(\mathbb{F})$}

Let $\Delta = {^2}A_{2n}(\mathbb{F})$, with $\mathbb{F}$ as in the previous subsection. In this case $U'_{A^*}$ is isomorphic to the additive group of $V(n^2,\mathbb{F})$ while $U_{A^*}/U'_{A^*}$ is isomorphic to the additive group of $V(n,\mathbb{K})$. The group $U'_{A^*}$ does not define a quotient of $\varepsilon_\Delta^n$. However, the function that maps every element $X\in\Gamma$ onto $(\varepsilon_\Delta^n(X)U'_{A^*})/U'_{A^*}$ is a projective embedding of $\Gamma$ in $U_{A^*}/U'_{A^*} \cong V(n,\mathbb{K})$, isomorphic to the natural embedding of $\Gamma$ in $V(n,\mathbb{K})$.

\subsubsection{The case $\Delta = {^2}D_{n+1}(\mathbb{K})$}

Let $\Delta = {^2}D_{n+1}(\mathbb{K})$. The discussion of Subsection \ref{BnK} can be repeated for this case with only a few minor changes. The elements of $U_{A^*}$ can be regarded as pairs $(S+M^TM, M)$, where $S$ is an $n\times n$ skew-symmetric matrix over $\mathbb{K}$ and $M$ a $2\times n$ matrix over $\mathbb{K}$. The pairs $(S,O)$ (where $O$ is the null $2\times n$ matrix) form a subgroup $W$ stabilized by $G_{A^*}$ and $U_{A^*}/W$ is a copy of $V\oplus V$, where $V:= V(n,\mathbb{K})$. Condition $(\mathrm{Q}1)$ of Subsection \ref{morphisms quotients} fails to hold for $W$. Hence $W$ does not define a quotient of $\varepsilon^n_\Delta$. Nevertheless, if we map every element $X\in\Gamma$ onto $(\varepsilon_\Delta^n(X)W)/W$ then we obtain a projective embedding of $\Gamma$ in $U_{A^*}/W\cong V\oplus V$. 

If $\mathrm{char}(\mathbb{K})\neq 2$ then $U'_{A^*} = W$. Let $\mathrm{char}(\mathbb{K}) = 2$. Then $U_{A^*} = 1$ and we can put a $\mathbb{K}$-vector space structure $V'$ on $U_{A^*}$ with scalar multiplication defined as follows: $t\cdot (S+M^TM,M) := (t^2(S+M^TM), tM)$. Denoted by $\mathbb{F}$ the subfield of square elements of $\mathbb{K}$, we have
\[\mathrm{dim}(V') = 2n + {{n-1}\choose 2}\cdot[\mathbb{K}:\mathbb{F}].\]
The embedding $\varepsilon_\Delta^n$ is $\mathbb{K}$-linear but non-projective. We have $\mathrm{dim}(\varepsilon_\Delta^n(p)) = 2$ for every point $p$ of $\Gamma$ while if $L$ is a line of $\Gamma$ then 
$\mathrm{dim}(\varepsilon_\Delta^n(L)) = 4 + [\mathbb{K}:\mathbb{F}]$. 

\subsubsection{Dominancy of $\varepsilon_\Delta^n$}\label{Back to general}

Turning back to the general case, let $\Delta$ be an arbitrary thick building of type $C_n$ with $n\geq 3$ and $A^*$ an $n$-elements of $\Delta$. The following is known (see e.g. \cite{Pas-Far}):

\begin{lemma}\label{Far polar lemma}
The geometry $\mathrm{Far}_\Delta(A^*)$ is simply connected except when $\Delta = C_3(2)$ $(= B_3(2)$) or
$\Delta = {^2}A_5(2)$. The universal cover of $\mathrm{Far}_\Delta(A^*)$ is a double cover when $\Delta = C_3(2)$ and a $4$-fold cover when $\Delta = {^2}A_5(2)$.
\end{lemma}

\begin{co}\label{Far polar cor}
The geometry $\mathrm{Far}_\Delta(A^*)$ is $2$-simply connected except when $\Delta$ is $C_3(2)$ or ${^2}A_5(2)$ and possibly when $\Delta$ is $C_n(2)$ or ${^2}A_{2n-1}(2)$, with $n > 3$.  
\end{co}

Corollary \ref{Far polar cor} immediately follows from Lemma \ref{Far polar lemma}. By Corollary \ref{Far polar cor} and Propositions \ref{truncation} and \ref{expansion} we obtain the following:

\begin{prop}\label{Unipotent polar}
The embedding $\varepsilon^n_\Delta$ is dominant except when $\Delta = C_3(2)$ or $\Delta = {^2}A_5(2)$ and possibly when $\Delta$ is $C_n(2)$ or ${^2}A_{2n-1}(2)$ with $n > 3$. 
\end{prop}

It is proved in Baumeister, Meixner and Pasini \cite{BMP} that when $\Delta = C_n(2)$ the codomain $U(\varepsilon_\Delta^n)$ of $\varepsilon_\Delta^n$ is elementary abelian of order $2^{2^n-1}$ ($> {{n+1}\choose 2} = |U_{A^*}|$). Hence in this case $\varepsilon_\Delta^n$ is not dominant, for any $n \geq 3$. Accordingly, the universal $2$-cover of $\mathrm{Far}_\Delta(A^*)$ is a $2^{2^n-1}/{{n+1}\choose 2}$-fold $2$-cover.

Let $\Delta = {^2}A_5(2)$. As stated in Lemma \ref{Far polar lemma}, the geometry $\mathrm{Far}_\Delta(A^*)$ admits a $4$-fold cover. Accordingly $|U(\varepsilon_\Delta^3)| = 2^{11} = 4\cdot|U_{A^*}|$. It is quite natural to conjecture that when $\Delta = {^2}A_{2n-1}(2)$ with $n > 3$ the embedding $\varepsilon_\Delta^n$ is not dominant. Equivalently,  $\mathrm{Far}_\Delta(A^*)_n$ is not $2$-simply connected.  

\section{A selection of results}\label{Results}

Throughout this section $\Gamma = ({\cal P},{\cal L})$ is (the $\{1,2\}$-truncation of) a shadow geometry of a thick building $\Delta^\circ$ of rank $n-1\geq 2$ and $\varepsilon:\Gamma\rightarrow V$ is a projective embedding of $\Gamma$. Moreover $\tilde{\varepsilon}$ is the hull of $\varepsilon$ (as in Subsection \ref{Hull}) and $\varepsilon^k_\Delta$ is the unipotent representation of $\Gamma$  of type $(\Delta, k)$ for a suitable building $\Delta$ of rank $n$ and a type $k$ of $\Delta$. In each of the cases to be examined in the sequel both conditions $(*)$ of Proposition \ref{connessione} and $(\mathrm{A})$ of Subsection \ref{sezione-epsilonA} hold, hence $\varepsilon^k_\Delta$ is well defined (by $(\mathrm{A})$) and it is an embedding (by $(*)$). In any case the embedding $\varepsilon^k_\Delta$ admits abelian quotients, hence we can consider its abelianization $\varepsilon^k_{\Delta,\mathrm{ab}}$. 

We firstly state a number of theorems proved in \cite{PasEE}. Next we state two new theorems, to be proved in Sections \ref{Section3} and \ref{Section4}. Finally, we prove a few corollaries on certain far-away geometries. 

We refer to Section \ref{notation} for the way of labelliing the nodes of the diagrams and for the names of particular buildings and their shadow geometries. 

\begin{theo}[\rm{\cite[Section 6]{PasEE}}]\label{ThA1} 
Let $\Delta^\circ$ be of type $C_{n-1}$ or $D_{n-1}$ and $\Gamma = \Delta^\circ_1$ (namely, $\Gamma$ is a polar space). Let $\Delta$ be the (uniquely determined) building of type $C_n$ or $D_n$ respectively such that $\Delta^\circ\cong \mathrm{Res}_\Delta(a)$ for a $1$-element $a$ of $\Delta$. Then all the following hold:\\
(1) The hull $\tilde{\varepsilon}$ of $\varepsilon$ is isomorphic to  $\varepsilon_\Delta^1$. \\
(2) The abelianization $\varepsilon^1_{\Delta,\mathrm{ab}}$ of $\varepsilon^1_\Delta$ is the absolute projective embedding of $\Gamma$.\\
(3) We have $\varepsilon = \tilde{\varepsilon}$ if and only if $\varepsilon$ embeds $\Gamma$ as a quadric in $\mathrm{PG}(V)$.  
\end{theo} 

\begin{theo}[\rm{\cite[Section 8.1]{PasEE}}]\label{ThA2} 
Given a field $\mathbb{F}$, let $\Gamma = A_{n-1,2}(\mathbb{F})$ (the line-grassmannian of $\mathrm{PG}(n-1,\mathbb{F})$). Let $\varepsilon$ be the usual embedding of $\Gamma$ into the exterior square $V = \wedge^2V_0$ of $V_0 := V(n,\mathbb{F})$, mapping every $2$-subspace $\langle v_1, v_2\rangle$ of $V_0$ onto the $1$-subspace $\langle v_1\wedge v_2\rangle$ of $V$. Let $\Delta = D_n(\mathbb{F})$. Then $\varepsilon = \tilde{\varepsilon}\cong \varepsilon_\Delta^n$.
\end{theo} 

\begin{theo}[\rm{\cite[Section 10]{PasEE}}]\label{ThA3}
Let $\Gamma = D_{5,5}(\mathbb{F})$ (the half-spin geometry of $\Delta^\circ = D_5(\mathbb{F})$) and let $\varepsilon$ be the half-spin embedding of $\Gamma$ into $V(16,\mathbb{F})$. Let $\Delta = E_6(\mathbb{F})$. Then $\varepsilon = \tilde{\varepsilon}\cong \varepsilon_\Delta^1$. 
\end{theo} 

\begin{theo}[\rm{\cite[Section 9.4]{PasEE}}]\label{ThA4}
Let $\Gamma = C_{3,3}(p)$ for a prime $p$ and $\Delta = F_4(p)$. Let $\varepsilon$ be the usual projective embedding of $\Gamma$ in a $14$-dimensional subspace of the third exterior power $\wedge^3V_0$ of $V_0 := V(6,p)$, induced by the natural embedding of $A_{5,3}(p)$ in $\mathrm{PG}(\wedge^3V_0)$. Then $\varepsilon = \varepsilon^4_{\Delta,\mathrm{ab}}$. Moreover:\\
(1) Let $p > 2$. Then $\tilde{\varepsilon} \cong \varepsilon_\Delta^4$.\\
(2) Let $p = 2$. Then $\varepsilon_\Delta^4$ is a quotient of $\tilde{\varepsilon}$ by a group of order $2$. 
\end{theo} 

\begin{theo}[\rm{\cite[Section 9.4]{PasEE}}]\label{ThA5}
Let $\Gamma = {^2}A_{5,3}(p)$ for a prime $p$ and $\Delta = {^2}E_6(p)$. Let $\varepsilon$ be the usual projective embedding of $\Gamma$ in $V = V(20,p)$, induced by the natural embedding of $A_{5,3}(p^2)$ in $\mathrm{PG}(\wedge^3V_0)$, $V_0 = V(5,p^2)$ (see e.g. Cooperstein and Shult {\rm \cite{CS-2001}}). Then $\varepsilon = \varepsilon^4_{\Delta,\mathrm{ab}}$. Moreover:\\
(1) Let $p > 2$. Then $\tilde{\varepsilon} \cong \varepsilon_\Delta^4$.\\
(2) Let $p = 2$. Then $\varepsilon_\Delta^4$ is a quotient of $\tilde{\varepsilon}$ by a group of order $4$. 
\end{theo}

We recall that $E_{6,1}(\mathbb{F})$ and $E_{7,7}(\mathbb{F})$ admit projective embeddings of dimension $27$ cand $56$ respectively. We refer to Cohen \cite[5.2]{Cohen-Handbook} (also Pasini \cite[5.1]{PasF4}) for a description of the $27$-dimensional projective embedding of $E_{6,1}(\mathbb{F})$ and to Cooperstein \cite{CoopE7} for the $56$-dimensional projective embedding of $E_{7,7}(\mathbb{F})$. The dimensions of these emebddings are equal to the size of a generating set of $E_{6,1}(\mathbb{F})$ and $E_{7,7}(\mathbb{F})$ respectively (Blok and Brouwer \cite{BB-span}, Cooperstein and Shult \cite{CS}).  Hence these two embeddings are linearly dominant. Being linearly dominant, they are absolute, since both $E_{6,1}(\mathbb{F})$ and $E_{7,7}(\mathbb{F})$ admit the absolute projective embedding (Kasikova and Shult \cite{KS}). Being absolute, they are uniquely determined by their dimensions.

The next two theorems will be proved in Sections \ref{Section3} and \ref{Section4} of this paper. 

\begin{theo}\label{ThB1}
Let $\Gamma = E_{6,1}(\mathbb{F})$ and let $\varepsilon$ be the $27$-dimensional projective embedding of $\Gamma$ mentioned above. Let $\Delta = E_7(\mathbb{F})$. Then $\varepsilon = \tilde{\varepsilon}\cong \varepsilon_\Delta^7$. 
\end{theo} 

\begin{theo}\label{ThB2}
Let $\Gamma = E_{7,7}(p)$ for a prime $p$ and let $\varepsilon$ be the $56$-dimensional projective embedding of $\Gamma$ mentioned above. Let $\Delta = E_8(p)$. Then $\tilde{\varepsilon}\cong \varepsilon_\Delta^8$ and $\varepsilon \cong \varepsilon_{\Delta,\mathrm{ab}}^8$.  
\end{theo} 

The proofs of Theorems \ref{ThA1}, \ref{ThA2} and \ref{ThA3} given in \cite{PasEE} run as follows. Firstly, it is proved that $\mathrm{Ex(\varepsilon)}$ is a (possibly improper) quotient of $\mathrm{Tr}^{\{1,2,3\}}(\mathrm{Far}_\Delta(a^*)_k)$, for $a^*$ an element of type $k^*$ with $(k, k^*)$ equal to $(1,1)$ in Theorem \ref{ThA1}, equal to $(n,n)$ or $(n,n-1)$ in Theorem \ref{ThA2} and equal to $(1,6)$ in Theorem \ref{ThA3}. On the other hand, $\mathrm{Tr}^{\{1,2,3\}}(\mathrm{Far}_\Delta(a^*)_k)\cong \mathrm{Ex}(\varepsilon_\Delta^k)$ by Theorem \ref{1.5}. Moreover, in each of the cases considered in Theorems \ref{ThA1}, \ref{ThA2} and \ref{ThA3} the geometry $\mathrm{Tr}^{\{1,2,3\}}(\mathrm{Far}_\Delta(a^*)_k)$ is simply connected (Cuypers and Pasini \cite{CP} and Pasini \cite{PasFar}). The conclusions are drawn from this fact by exploiting Proposition \ref{expansion}. 

The proofs of Theorems \ref{ThA4} and \ref{ThA5} are quite different. We are not going to recall them here. We only remark that the simple connectedness of $\mathrm{Tr}^{\{1,2,3\}}(\mathrm{Far}_\Delta(a^*)_k)$ is not exploited in those proofs. On the contrary, it is obtained from those theorems. Explicitly, by Theorems \ref{ThA4} and \ref{ThA5}, Propositions \ref{expansion} and \ref{truncation} and Theorem \ref{1.5} we obtain the following:

\begin{co}
Let $\Delta = F_4(p)$ for a prime $p$ and let $a^*$ be a $4$-element of $\Delta$.\\
(1) Let $p > 2$. Then $\mathrm{Far}_\Delta(a^*)_4$ is $2$-simply connected. Its $\{1,2,3\}$-truncation $\mathrm{Tr}^{\{1,2,3\}}(\mathrm{Far}_\Delta(a^*)_4)$ is simply connected. \\
(2) Let $p = 2$. Then $\mathrm{Far}_\Delta(a^*)_4$ admits a $2$-fold $2$-cover. The $\{1,2,3\}$-truncation of that $2$-cover is a $2$-fold cover of $\mathrm{Tr}^{\{1,2,3\}}(\mathrm{Far}_\Delta(a^*)_4)$. 
\end{co} 

\begin{co}
Let $\Delta = {^2}E_6(p)$ for a prime $p$. Let $a^*$ be a $4$-element of $\Delta$.\\
(1) Let $p > 2$. Then $\mathrm{Far}_\Delta(a^*)_4$ is $2$-simply connected. Its $\{1,2,3\}$-truncation $\mathrm{Tr}^{\{1,2,3\}}(\mathrm{Far}_\Delta(a^*)_4)$ is simply connected. \\
(2) Let $p = 2$. Then $\mathrm{Far}_\Delta(a^*)_4$ admits a $4$-fold $2$-cover. The $\{1,2,3\}$-truncation of that $2$-cover is a $4$-fold cover of $\mathrm{Tr}^{\{1,2,3\}}(\mathrm{Far}_\Delta(a^*)_4)$. 
\end{co} 

In a similar way, we obtain the following from Theorems \ref{ThB1} and \ref{ThB2}.

\begin{co}
Let $\Delta = E_7(\mathbb{F})$ and let $a^*$ be a $7$-element of $\Delta$. Then $\mathrm{Far}_\Delta(a^*)_4$ is $2$-simply connected. Accordingly, $\mathrm{Tr}^{\{1,2,3\}}(\mathrm{Far}_\Delta(a^*)_4)$ is simply connected.
\end{co}

\begin{co}
Let $\Delta = E_8(p)$ for a prime $p$ and let $a^*$ be an $8$-element of $\Delta$. Then $\mathrm{Far}_\Delta(a^*)_4$ is $2$-simply connected. Accordingly, $T^{\{1,2,3\}}(\mathrm{Far}_\Delta(a^*)_4)$ is simply connected.
\end{co}

\section{Conjectures and problems} 

(1) Let $\varepsilon$ be the half-spin embedding of $D_{n,n}(\mathbb{F})$. We conjecture that when $n = 6$ or $7$ then $\tilde{\varepsilon}\cong \varepsilon_\Delta^1$ for $\Delta = E_{n+1}(\mathbb{F})$. When $n > 7$ then $\mathrm{Ex}(\varepsilon)$ is a truncation of a geometry belonging to the following non-spherical-like diagram of rank $n+1 > 8$:

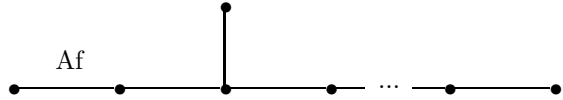
\begin{picture}(310,53)(0,0)
\put(50,8){$\bullet$}
\put(53,11){\line(1,0){37}}
\put(68,18){$\mathrm{Af}$}
\put(90,8){$\bullet$}
\put(93,11){\line(1,0){37}}
\put(130,8){$\bullet$}
\put(133,11){\line(1,0){37}}
\put(170,8){$\bullet$}
\put(173,11){\line(1,0){12}}
\put(190,11){$...$}
\put(203,11){\line(1,0){12}}
\put(215,8){$\bullet$}
\put(218,11){\line(1,0){37}}
\put(255,8){$\bullet$}
\put(132,10){\line(0,1){30}}
\put(130,39){$\bullet$}
\end{picture}

\noindent
Most likely, the expansion $\mathrm{Ex}(\tilde{\varepsilon})$ of the hull $\tilde{\varepsilon}$ of $\varepsilon$ (which, by Proposition \ref{expansion}, is isomorphic to the universal cover of $\mathrm{Ex}(\varepsilon)$) has infinite diameter.

Situations like this occur in many other cases, as when $\Gamma$ is of type $C_{n,n}$ with $n > 3$ or of type $F_{4, k}$ with $k = 1$ or $4$, or $E_{n, k}$ with either $k = 2$ or $(n, k) = (7,1), (8,1)$ or $(8,8)$. 

\bigskip

\noindent
(2) Perhaps the conclusions of Theorems \ref{ThA4}, \ref{ThA5} and \ref{ThB2} hold in general, with the prime field $\mathbb{F}_p$ replaced any field $\mathbb{F}$ ($\neq \mathbb{F}_2$ in Theorems \ref{ThA4} and \ref{ThA5}).  

\bigskip

\noindent
(3) As for dual polar spaces of rank $3$,  we have only considered $C_{3,3}(\mathbb{F})$ and ${^2}A_{5,3}(\mathbb{F})$, moreover with $\mathbb{F} = \mathbb{F}_p$.  What can we say on embeddable dual polar spaces of types different from these? 

For instance, let $\varepsilon$ be the spin embedding of $\Gamma = B_{3,3}(\mathbb{F})$ in $V(8,\mathbb{F})$.With the help of a result by Cuypers and Van Bon \cite{vBC}, it is proved in \cite[Section 9.5]{PasEE} that $\varepsilon$ is a quotient of $\varepsilon_\Delta^1$ with $\Delta = F_4(\mathbb{F})$. When $\mathrm{char}(\mathbb{F})\neq 2$ then $\varepsilon\cong \varepsilon_{\Delta,\mathrm{ab}}^1$. When $\mathrm{char}(\mathbb{F}) = 2$ then $\varepsilon_{\Delta,\mathrm{ab}}^1$ is linear but it is larger than $\varepsilon$. Its codomain is (the additive group of) a $14$-dimensional module for $\mathrm{Spin}(7,\mathbb{F})$. The spin module is a factor of that module. 

With $\varepsilon$, $\Gamma$ and $\Delta$ as in the previous paragraph, if $\mathbb{F} = \mathbb{F}_2$ then $\varepsilon^1_\Delta$ is not dominant, as it follows from Theorem \ref{ThA4} recalling that $B_{3,3}(2)\cong C_{3,3}(2)$. We conjecture that if $\mathbb{F}\neq \mathbb{F}_2$ then $\varepsilon_\Delta^1$ is dominant.  

\bigskip

\noindent
(4) Let $\mathbb{D}$ be a Cayley division algebra over a field $\mathbb{F}$ and let $\Pi$ be the polar space of rank $3$ with planes defined over $\mathbb{D}$ (Tits \cite[Chapter 9]{Tits}). It is well known that $\Pi$ can be constructed as a subgeometry of $E_7(\mathbb{F})$, the points, lines and planes of $\Pi$ being certain elements of $E_7(\mathbb{F})$ of type $1$, $6$ and $7$ respectively. 

The polar space $\Pi$ does not admit any projective embedding. However, the dual $\Gamma$ of $\Pi$ admits a $56$-dimensional projective embedding $\varepsilon_\Gamma$, induced by the projective embedding of $E_{7,7}(\mathbb{F})$ in $V(56,\mathbb{F})$. We refer to De Bruyn and Van Maldeghem \cite{DBVM} for a thorough investigation of $\varepsilon_\Gamma$. In particular, $\varepsilon_\Gamma$ is linearly dominant (hence absolute, in view of \cite{KS}). We wonder what the hull of $\varepsilon_\Gamma$ might be. 

\section{Generators and more on embeddings}\label{preliminaries}

In this section we collect a few definitions and general results to be exploited later in Sections \ref{Section3} and \ref{Section4}, in the proofs of Theorems \ref{ThB1} and \ref{ThB2}.

\subsection{Subspaces, generators and hyperplanes}\label{generating}

Let $\Gamma = ({\cal P},{\cal L})$ be a point-line geometry. A subset $S\subseteq{\cal P}$ is called a {\em subspace} if it contains every line $L\in{\cal L}$ such that $|L\cap S|> 1$. Obviously, $\cal P$ is a subspace of $\Gamma$. A subspace is said to be {\em proper} if it is different from ${\cal P}$. 

The intersection of a family of susbpaces of $\Gamma$ is a subspace. Given a subset $S\subseteq {\cal P}$, the subspace $\langle S\rangle_\Gamma$ {\em generated} by $S$ is the smallest susbspace of $\Gamma$ containing $S$. It is the intersection of all subspaces of $\Gamma$ that contain $S$. It can also be described as follows. Put $S_0 := S$ and let $S_{n+1}$ be the union $S_n$ and all lines of $\Gamma$ that meet $S_n$ in at least two points. Then 
\begin{equation}\label{generating-equation}
\langle S\rangle_\Gamma = \bigcup_{n=0}^\infty S_n.
\end{equation}
A subset $S\subseteq {\cal P}$ is said to {\em generate} $\Gamma$ if $\langle S\rangle_\Gamma = {\cal P}$. If $S$ generates $\Gamma$ then we also say that $S$ is {\em generating set} of $\Gamma$. 

A {\em hyperplane} of $\Gamma$ is a proper subspace $H$ of $\Gamma$ such that every line of $\Gamma$ meets $H$ non-trivially. It is easy to see that a hyperplane $H$ is a maximal subspace of $\Gamma$ if and only if the collinearity graph of $\Gamma$ induces a connected graph on ${\cal P}\setminus H$. 

\subsection{Generating sets and embeddings}\label{generating-embedding}

Let $\varepsilon:\Gamma\rightarrow G$ be an embedding. The following is a straightforward consequence of equality (\ref{generating-equation}) of Subsection \ref{generating}. 

\begin{lemma}\label{generating-lemma}
Let $S$ be a generating set of $\Gamma$. Then $G = \langle \varepsilon(p)\rangle_{p\in S}$.
\end{lemma}

\begin{prop}\label{generating-corollary}
Suppose that $\varepsilon$ as well as its hull $\tilde{\varepsilon}$ are abelian. Suppose moreover that $G = \oplus_{p\in S}\varepsilon(p)$ (direct sum of abelian groups) for a generating set $S$ of $\Gamma$. Then $\tilde{\varepsilon} = \varepsilon$.  
\end{prop}
\pr By Lemma \ref{generating-lemma}, the codomain $U(\varepsilon)$ of $\tilde{\varepsilon}$ is generated by the subgroups $\tilde{\varepsilon}(p)$ for $p\in S$. By assumption, $U(\varepsilon)$ is commutative and $G = \oplus_{p\in S}\varepsilon(p)$. It follows that the projection $\pi_\varepsilon:U(\varepsilon)\rightarrow G$ is an isomorphism. \eop

\subsection{More on absolute projective embeddings}

Suppose that $\Gamma$ admits the absolute projective embedding and let $\varepsilon:\Gamma\rightarrow V$ be its absolute projective embedding. It readily follows from the definition of absolute projective embeddings that $\varepsilon$ is homogeneous. Denoted by $G$ the automorphism group of $\Gamma$, let $\widehat{G}$ be the lifting of $G$ to $V$ through $\varepsilon$. Let $U$ be a vector subspace of $V$ defining a quotient of $\varepsilon$. As noticed in Subsection \ref{Homogeneity}, if $U$ is stabilized by $\widehat{G}$ then the embedding $\varepsilon/U$ is homogeneous. Since $\varepsilon$ is abolute, the converse also holds true (Pasini and Van Maldeghem \cite[Proposition 13]{PVM}):

\begin{prop}\label{homogeneous quotient absolute}
If $\varepsilon/U$ is homogeneous then $\widehat{G}$ stabilizes $U$.
\end{prop}

\section{Proof of Theorem \ref{ThB1}}\label{Section3}

Let $\Delta^\circ = E_6(\mathbb{F})$ and $\Gamma = E_{6,1}(\mathbb{F})$ for a field $\mathbb{F}$ and let $\varepsilon:\Gamma\rightarrow V := V(27,\mathbb{F})$ be the $27$-dimensional projective embedding of $\Gamma$. We recall that $\varepsilon$ is absolute. Let $\tilde{\varepsilon}$ be the hull of $\varepsilon$.

\begin{lemma}\label{E6-abelian}
The embedding $\tilde{\varepsilon}$ is abelian.
\end{lemma}
\pr Let $S$ be a $6$-element of $\Delta^\circ$ and $P(S)$ the set of points of $\Gamma$ incident to $S$. The set $P(S)$ is a subspace of $\Gamma$ isomorphic $D_{5,1}(\mathbb{F})$. The embeddings $\tilde{\varepsilon}$ and $\varepsilon$ induce on $P(S)$ embeddings $\tilde{\varepsilon}_S$ and $\varepsilon_S$ such that $\varepsilon_S$ is a quotient of $\tilde{\varepsilon}_S$ and $\varepsilon_S$ is isomorphic to the natural emebdding of $D_{5,1}(\mathbb{F})$ in $V(10,\mathbb{F})$. By claim (3) of Theorem \ref{ThA1}, $\tilde{\varepsilon}_S = \varepsilon_S$. In particular, $\tilde{\varepsilon}_S$ is abelian. 

It is well known that any two points $x$ and $y$ of $\Gamma$ are incident to a common $6$-element of $\Delta^\circ$ (see e.g. Cohen and Cooperstein \cite{CC}, also \cite[7.6.1]{PasDG}). Therefore $[\tilde{\varepsilon}(x),\tilde{\varepsilon}(y)] = 1$ by the previous paragraph. Hence $\tilde{\varepsilon}$ is abelian. \eop 

\begin{prop}\label{E6-dominant}
We have $\tilde{\varepsilon} = \varepsilon$.
\end{prop}
\pr The geometry $\Gamma$ admits a generating set $S$ of size $27$, formed by the points of $\Gamma$ contained in a given apartment of $\Delta^\circ$ (Blok and Brouwer \cite{BB-span}, Cooperstein and Shult \cite{CS}). Hence $V = \oplus_{p\in S}\varepsilon(p)$, by Lemma \ref{generating-lemma} and since $\mathrm{dim}(V) = 27$. The conclusion follows from Lemma \ref{E6-abelian} and Proposition \ref{generating-corollary}. \eop

\bigskip

The automorphism group $G^\circ := \mathrm{Aut}(\Gamma)$ of $\Gamma$ is isomorphic to the group $\mathrm{Aut}(\Delta^\circ)$ of type-preserving automorphisms of $\Delta^\circ$. Let $\widehat{G}^\circ$ be the lifting of $G^\circ$ to $V$ through $\varepsilon$ (recall that $\varepsilon$, being absolute, is homogeneous). Obviously, $\widehat{G}^\circ = Z\cdot G^\circ$, where $Z$ is the group of scalar transformations of $V$. Clearly, $V$ is a $\widehat{G}^\circ$-module.

The next proposition is not necessary for proof of Theorem \ref{ThB1}, but we will use it in Section \ref{Section4}, in the proof of Theorem \ref{ThB2}. Moreover, it is interesting in itself. 

\begin{prop}\label{E6-irreducibility}
The $\widehat{G}^\circ$-module $V$ is irreducible.
\end{prop}
\pr A cubic form $f$ is defined on $V$ such that the set $S_0 := \{v\in V~|~ f(v) = 0\}$ contains $\cup_{X\in\Gamma}\varepsilon(X)$ and  two vectors $u, v \in V\setminus S_0$ belong to the same orbit of $\widehat{G}^\circ$ if and only if $f(v)^{-1}f(u)\in \mathbb{F}^{*(3)}$, where $\mathbb{F}^{*(3)} := \{t^3\}_{t\in \mathbb{F}^*}$ and $\mathbb{F}^*$ is the multiplicative group of $\mathbb{F}$ (Cohen \cite[5.2]{Cohen-Handbook}, also Pasini \cite[Section 5]{PasF4}). The set $S_0$ is itself the union of a number of orbits of $\widehat{G}^\circ$. It is not so difficult to see that each of the orbits not contained in $S_0$ spans $V$ while every non-trivial orbit contained in $S_0$ contains vectors $v$ and $w$ such that $v+w\not\in S_0$. It follows that no proper non-trivial subspace of $V$ can be obtained as the union of $\{0\}$ and some of the orbits of $\widehat{G}^\circ$. Hence no non-trivial proper subspace of $V$ is stabilized by $\widehat{G}$.  \eop

\begin{co}\label{E6-irreducibility Cor}
The embedding $\varepsilon$ does not admit any proper homogeneous projective quotient.
\end{co}
\pr Since $\varepsilon$ is absolute, if a subspace $U$ of $V$ defines a homogeneous quotient of $\varepsilon$ then $U$ is stabilized by $\widehat{G}^\circ$, by Proposition \ref{homogeneous quotient absolute}. Hence $U = 0$, by Proposition \ref{E6-irreducibility}.  \eop

\bigskip

Let now $\Delta = E_7(\mathbb{F})$ and $a$ and $a^*$ be opposite elements of $\Delta$ of type $7$. Hence $\Gamma\cong \Delta_7(a)$. Let $\varepsilon_{\Delta}^7$ be the unipotent representation of $\Gamma$ in $U_{a^*}$.

\begin{lemma}\label{E6 Uastar projective}
The embedding $\varepsilon_\Delta^7$ is projective (over $\mathbb{F}$) and $27$-dimensional.  
\end{lemma}
\pr It is proved in Cooperstein \cite[(3.12)]{CoopE7} that $U_{a^*}$ is isomorphic to the additive group of $V(27,\mathbb{F})$. (We warn that Cooperstein assumes $\mathrm{char}(\mathbb{F})\neq 2$ in \cite{CoopE7}, but his description of $U_{a^*}$ remains valid when $\mathrm{char}(\mathbb{F}) = 2$.) Using the information offered in \cite{CoopE7}, it is not so difficult to see that, regarded $U_{a^*}$ as an $\mathbb{F}$-vector space, for every element $X\in \Gamma$ the subgroup $\varepsilon_\Delta^7(X)$ is a vector subspace of $U_{a^*}$, of dimension $1$ or $2$ according to whether $X$ is a point or a line of $\Gamma$. The lemma is proved.  \eop

\bigskip

\noindent
{\bf Remark.} We have based our proof of Lemma \ref{E6 Uastar projective} on \cite{CoopE7}, but a part of the statement of that lemma can be obtained in a more elementary way. Recall that any two points of $\Gamma$ are incident with a common $6$-element of $\Delta^\circ$. Let $S$ be a $6$-element of $\Delta^\circ$ and $P(S)$ the set of points of $\Gamma$ incident to $S$. Let $\varepsilon_{\Delta,S}^7$ be the embedding induced by $\varepsilon_\Delta^7$ on $P(S)$. By Theorem \ref{unipotente indotto} the induced embedding $\varepsilon_{\Delta,S}^7$ is isomorphic to the unipotent representation of $D_{5,1}(\mathbb{F})$ relative to the pair $(D_6(\mathbb{F}), 1)$. By Claim (3) of Theorem \ref{ThA1} the latter embedding is projective, whence abelian. Recall now that any two points of $\Gamma$ are incident with a common $6$-element of $\Delta^\circ$. It follows that $\varepsilon_\Delta^7$ is locally projective and abelian. 

\begin{prop}\label{EndThB1}
We have $\varepsilon_\Delta^7 \cong \varepsilon$.
\end{prop}
\pr Since $\varepsilon$ is absolute and $\varepsilon_\Delta^7$ is projective, $\varepsilon_\Delta^7$ is a quotient of $\varepsilon$. The isomorphism $\varepsilon_\Delta^7\cong \varepsilon$ follows from the fact that $\varepsilon_\Delta^7$ has the same dimension as $\varepsilon$. Alternatively, we can remark that $\varepsilon_\Delta^7$ is homogeneous by Theorem \ref{epsilonA-omogeneo} and apply Proposition \ref{E6-irreducibility}.  \eop 

\bigskip

Propositions \ref{E6-dominant} and \ref{EndThB1} give us the statement of Theorem \ref{ThB1}.  

\section{Proof of Theorem \ref{ThB2}}\label{Section4}

Given a field $\mathbb{F}$, let $\Delta^\circ = E_7(\mathbb{F})$, $\Gamma = E_{7,7}(\mathbb{F})$ and $G^\circ := \mathrm{Aut}(\Gamma)$ ($\cong \mathrm{Aut}(\Delta^\circ)$). Let $V = V(27,\mathbb{F})$ and  let $\varepsilon:\Gamma\rightarrow V$ be projective embedding of $\Gamma$ in $V$. As remarked in Section 4, the embedding $\varepsilon$ is absolute. Hence it is homogeneous. We denote by $\widehat{G}^\circ$ the lifting of $G^\circ$ to $V$ through $\varepsilon$. So, the vector space $V$ is in fact a $\widehat{G}^\circ$-module.  

We denote the point-set and the line-set of $\Gamma$ by $\cal P$ and $\cal L$, as usual. We also denote by $d$ the distance in the collinearity graph of $\Gamma$ and for two points $x, y\in {\cal P}$ we write $x\perp y$ when $d(x,y) \leq 1$. As usual, $x^\perp$ stands for the set of points at distance at most $1$ from $x$.

We recall that the collinearity graph of $\Gamma$ has diameter equal to $3$. Moreover, if $d(x,y) = 2$ for two points $x, y\in {\cal P}$, then there exists a $1$-element $S$ of $\Delta^\circ$ such that $x, y\in P(S) := \mathrm{sh}_7(S)$ (see Cohen and Cooperstein \cite{CC}). 

\subsection{A property of the $\widehat{G}^\circ$-module $V$} 

In this subsection we state a result on the $\widehat{G}^\circ$-module $V$, valid for any choice of $\mathbb{F}$. We shall turn to the special case $\mathbb{F} = \mathbb{F}_p$ in the next subsection, where we shall prove Theorem \ref{ThB2}.    

\begin{prop}\label{E7 irreducible}
The $\widehat{G}^\circ$-module $V$ is irreducible. 
\end{prop}
\pr Let $U$ be a subspace of $V$ stabilized by $\widehat{G}^\circ$. Given two points $a_1$ and $a_2$ of $\Gamma$ at distance $3$, for $i = 1, 2$ let $V(a_i) := \langle \varepsilon(x)\rangle_{x\in a_i^\perp}$. Then $\mathrm{dim}(V(a_1)) = \mathrm{dim}(V(a_2)) = 28$, we have $V = V(a_1)\oplus V(a_2)$ and for $i = 1, 2$ the embedding $\varepsilon$ induces on $\overline{V}_i := V(a_i)/\varepsilon(a_i)$ a projective embedding $\varepsilon_i$ isomorphic to the $27$-dimensional projective embedding of $\mathrm{Res}_{\Delta^\circ}(a_i)\cong E_{6,1}(\mathbb{F})$ (see Cooperstein \cite{CoopE7}). The stabilizer $\widehat{G}^\circ_{a_1,a_2}$ of $\varepsilon(a_1)$ and $\varepsilon(a_2)$ induces on $\overline{V}_i$ a lifting $\overline{G}_i$ of $\mathrm{Aut}(E_{6,1}(\mathbb{F})$.   

As $\widehat{G}^\circ$ stabilizes $U$, for $i = 1, 2$ the group $\widehat{G}_{a_1,a_2}$ stabilizes the projection $U_i$ of $U$ onto $V(a_i)$ induced by the projection of $V$ onto $V(a_i)$. It follows that $\overline{U}_i := (U_i\varepsilon(a_i))/\varepsilon(a_i)$ is stabilized by $\overline{G}_i$. By Proposition \ref{E6-irreducibility} either $U_i\leq \varepsilon(a_i)$ or $U_i\varepsilon(a_i) = V(a_i)$. However this must be the case for any choice of two points $a_1$ and $a_2$ of $\Gamma$ at distance $3$. It is not so difficult to check that this can happen only if either $U = 0$ or $U = V$.  \eop

\begin{co}\label{0}
The embedding $\varepsilon$ does not admit any homogeneous projective proper quotient.
\end{co}
\pr Easy, by Propositions \ref{homogeneous quotient absolute} and \ref{E7 irreducible}.  \eop 

\subsection{The case $\mathbb{F} = \mathbb{F}_p$}

Henceforth $\mathbb{F} = \mathbb{F}_p$ for a prime $p$. Thus, $\Delta^\circ = E_7(p)$, $\Gamma = E_{7,7}(p)$ and $V = V(56,p)$. Throughout this subsection $\eta:\Gamma\rightarrow W$ is a locally projective embedding of $\Gamma$ satifying both the following:

\begin{itemize}
\item[$(\eta 1)$] The embedding $\eta$ is homogeneous. 
\item[$(\eta 2)$] The group $W$ is non-perfect.
\end{itemize}
As $\eta$ is locally projective, $|\eta(x)| = p$ for every point $x\in{\cal P}$ and $\eta(X)$ is elementary abelian of order $p^2$ for every line $X\in{\cal L}$. Consequently, $W$ is generated by elements of order $p$ and the factor group $W/W'$ of $W$ over its derived group $W'$ is an elementary abelian $p$-group.

Recall that, as stated before, $G^\circ$ stands from $\mathrm{Aut}(\Gamma)$. As $\eta$ is assumed to be homogeneous, $G^\circ$ lifts through $\eta$ to a subgroup of $\mathrm{Aut}(W)$, henceforth denoted by $\widehat{G}^\circ$. 

Henceforth, for an element $X\in \Gamma$ we write $V_X$ instead of $\eta(X)$. 

\begin{lemma}\label{primo}
We have $W'\cap V_X = 1$ for every element $X\in\Gamma$.
\end{lemma}
\pr Suppose to the contrary that $W'\cap V_X \neq 1$ for some element $X\in \Gamma$. If $X$ is a point then $V_X\leq W'$, since $|V_X| = p$ and $p$ is prime. If $X$ is a line then every element of $V_X$ belongs to $V_x$ for some point $x\in X$. Thus, we can assume that $W'\geq V_x$ for some point $x$. The automorphism group $G^\circ$ of $\Gamma$ acts transitively on $\cal P$. Moreover, $\eta$ is homogeneous by assumption and $W'$ is a characteristic subgroup of $W$. It follows that $W'$ contains $V_x$ for every $x\in{\cal P}$. This forces $W' = W$, contrary to $(\eta 2)$. Therefore $W'\cap V_X = 1$ for every $X\in\Gamma$. \eop

\begin{lemma}\label{W/W'}
The derived subgroup $W'$ of $W$ defines a quotient of $\eta$ and we have $\eta/W'\cong \varepsilon$. 
\end{lemma}
\pr  By Lemma \ref{primo}, we have $W'\cap V_X = 1$ for every $X\in \Gamma$, namely $W'$ satisfies condition $(\mathrm{Q}1)$ of Subsection \ref{morphisms quotients}. We shall prove that it also satisfies condition $(\mathrm{Q}2)$, thus proving that it defines a quotient of $\eta$. 

By way of contradiction, let $W'V_x\cap W'V_y > W'$ for two distinct points $x, y\in{\cal P}$. Then there exist elements $w_x, w_y\in W'$, $v_x\in V_x\setminus\{1\}$ and $v_y\in V_y\setminus\{1\}$ such that $w_xv_x = w_yv_y$. Hence $w_y^{-1}w_x = v_yv_x^{-1}$. Moreover $v_yv_x^{-1} \neq 1$ because $V_x\cap V_y = 1$ (indeed $x\neq y$ by assumption) and $v_x\neq 1 \neq v_y$. 

Put $d_0 = d(x,y)$. Suppose that $d_0 = 1$ and let $X$ be the line of $\Gamma$ through $x$ and $y$. Then by $1 \neq w_y^{-1}w_x = v_yv_x^{-1}\in V_X$ we obtain that $W'\cap V_X \neq 1$, contrary to Lemma \ref{primo}. Therefore $d_0 \geq 2$. 

The equality $w_xv_x = w_yv_y$ implies that $v_y\in W'V_x$. However $v_y$ generates $V_y$ since $v_y \neq 1$ and $V_y$ is cyclic of prime order $p$. Therefeore $V_y\subset W'V_x$. (Recall that $W'V_y$ is a group, since $W'$ is normal in $W$.) 

For every $d \in \{0,1,2,3\}$, the group $G^\circ$ acts transitively on the set of ordered pairs of points of $\Gamma$ at distance $d$. Moreover its lifting $\widehat{G}^\circ$ stabilizes $W'$, since $W'$ is a characteristic subgroup of $W$. Therefore, considering the stabilizer of $x$ in $G$, we obtain that $W'V_x$ contains $V_y$ for every point $y$ at distance $d_0$ from $y$. 

We have previously proved that  $2 \leq d_0 \leq 3$. Assume firstly that $d_0 = 3$. The set of points of $\Gamma$ at distance $3$ from a given point generates $\Gamma$. Hence $W = \langle V_y\rangle_{d(y,x) = d_0}$ by Lemma \ref{generating-lemma}. It follows that $W'V_x = W$. In particular, $V_y\leq W'V_x$ for a point $y$ collinear with $x$, contrary to what we have previously proved. 

Therefore $d_0 = 2$. The set of points at distance $2$ from $x$ generates a hyperplane $H_x$ of $\Gamma$, formed by all points of $\Gamma$ at distance at most $2$ from $x$ (Blok and Brouwer \cite{BB}, also Cohen and Cooperstein \cite{CC}). Hence $W'V_x$ contains $V_y$ for every $y\in H$. However $H$ also contains all lines through $x$. Therefore we also have $W'V_x \geq V_y$ for every point $y\in x^\perp$. We have previously proved that this cannot be. 

Consequently $W'V_x\cap W'V_y = W'$ for any two distict points $x, y \in {\cal P}$. Thus, $W'$ defines a quotient of $\eta$. Since $\widehat{G}^\circ$ stabilizes $W'$, the quotient $\eta/W'$ is homogeneous. Moreover the factor group $W/W'$, being elementary abelian of exponent $p$, can be regarded as the additive group of an $\mathbb{F}_p$-vector space. Therefore the embedding $\eta/W'$ is projective. Hence it is a quotient of $\varepsilon$, as the latter is absolute. Corollary \ref{0} now implies that $\eta/W' \cong\varepsilon$. \eop 

\bigskip

By Lemma \ref{W/W'} we immediately obtain the following.

\begin{co}\label{W/W' cor}
The factor group $W/W'$ is elementary abelian of order $p^{56}$.
\end{co}

For every $1$-element $S$ of $\Delta^\circ$ the set $P(S)$ is a subspace of $\Gamma$ and the geometry $\Gamma_S$ induced by $\Gamma$ on $P(S)$ is isomorphic to $D_{6,1}(p)$. The mapping $\eta_S$ induced by $\eta$ on $\Gamma_S$ is an embedding of $\Gamma_S$ in the subgroup $W_S := \langle V_x\rangle_{x\in P(S)}$ of $W$.

\begin{lemma}\label{etaS}
For every $1$-element $S$ of $\Delta^\circ$, the embedding $\eta_S$ is isomorphic to the (unique) projective embedding of $\Gamma_S \cong D_{6,1}(p)$ in $V(12,p)$. 
\end{lemma}
\pr Put $W'_S :=  W_S\cap W'$. It follows from Lemma \ref{W/W'} that $W'_S$ defines a quotient $\eta_S/W'_S$ of $\eta_S$. The embedding $\eta_S/W'_S$ is projective. However $\Gamma_S \cong D_{6,1}(p)$ and $D_{6,1}(p)$ admits a unique projective embedding, which is $12$-dimensional. Hence $\eta_S$ is isomorphic to the unique projective embedding of $D_{6,1}(p)$ in $V(12,p)$.  \eop

\bigskip
  
The first claim of the next corollary immediately follows from the statement of Lemma \ref{etaS}. The second claim has been proved at the very end of the proof of Lemma \ref{etaS}.

\begin{co}\label{WS abelian}
For every $1$-element $S$ of $\Delta^\circ$ the group $W_S$ is elementary abelian of order $p^{12}$. Moreover $W_S\cap W' = 1$. 
\end{co}

For two points $x, y \in {\cal P}$ we denote by  $[V_x, V_y]$ the subgroup of $W'$ generated by the commutators 
$[u,v]$ for $u\in V_x$ and $v\in V_y$. 

\begin{co}\label{distance 2}
For two points $x, y\in{\cal P}$, if $d(x,y) \leq 2$ then $[V_x,V_y] = 1$.
\end{co}
\pr When $d(x,y) \leq 1$ the claim is obvious. Let $d(x,y) = 2$. Then there exists a $1$-element $S$ of $\Delta^\circ$ such that $x, y\in P(S)$. Accordingly, $V_x, V_y \subset W_S$. However $W_S$ is abelian, by Corollary \ref{WS abelian}. Hence $V_x$ and $V_y$ commute.  \eop 

\begin{lemma}\label{commutatori}
Let $x_1, x_2, y_1, y_2$ be points of $\Gamma$ with $d(x_1, y_1) = d(x_2, y_2) = 3$. Then $[V_{x_1},V_{y_1}] = [V_{x_2}, V_{y_2}]$.
\end{lemma}
\pr Let $x, y_1, y_2$ be points such that $d(x,y_1) = d(x,y_2) = 3$ but $d(y_1, y_2) = 1$. Let $X$ be the line of $\Gamma$ through $y_1$ and $y_2$. The set $H_x$ of points of $\Gamma$ at distance at most $2$ from $x$ is a hyperplane of $\Gamma$ (Block and Brouwer \cite{BB}, also Cohen and Cooperstein \cite{CC}). Hence $X\cap H_x$ contains a point, say $z$. For every choice of $b_1\in V_{y_1}\setminus\{1\}$ and $b_2\in V_{2}\setminus\{1\}$ there is an element $c\in V_z$ such that $b_2 = b_1c$. By Corollary \ref{distance 2}, $[a,c] = 1$ for every $a\in V_x$. Hence $[a,b_1] = [a,b_2]$. It follows that $[V_x,V_{y_1}] = [V_x,V_{y_2}]$.

By Proposition \ref{connessione}, the collinearity graph of $\Gamma$ induces a connected graph on ${\cal P}\setminus H_x$. Therefore $[V_x, V_{y_1}] = [V_x,V_{y_2}]$ for any two points $y_1$ and $y_2$ at distance $3$ from $x$. Put $C_x := [V_x,V_y]$, where $y$ is a point at distance $3$ from $x$, no matter which. 

Let now $x_1$ and $x_2$ be two collinear points. Let $X$ be the line through $x_1$ and $x_2$ and $z$ a point of $X\setminus\{x_1, x_2\}$. Then there exists a point $y$ at distance $2$ from $z$ and $3$ from either of $x_1$ and $x_2$ (see e.g. Cohen and Cooperstein \cite{CC}). As $d(y,x_1) = d(y,x_2) = 3$, we have $[V_y, V_{x_1}] = [V_y,V_{x_2}]$. However $[V_y,V_{x_i}] = [V_{x_i},V_y] = C_{x_i}$ for $i = 1, 2$. Hence $C_{x_1} = C_{x_2}$. 

By the connectedness of the collinearity graph of $\Gamma$, we have $C_{x_1} = C_{x_2}$ for any two points $x_1$ and $x_2$. Namely, $[V_{x_1},V_{y_1}] = [V_{x_2}, V_{y_2}]$ for any choice of points  $x_1, x_2, y_1, y_2$ with $d(x_1, y_1) = d(x_2, y_2) = 3$.  \eop 

\bigskip

Lemma \ref{commutatori} allows us to adopt the following notation: we put $C = [V_x,V_y]$ for a given pair of points at distance $3$, no matter which. Let $Z(W)$ be the center of $W$. 

\begin{co}\label{commutatori bis}
We have $W' = C \leq Z(W)$.
\end{co}
\pr The equality $W' = C$ is obvious. Let us prove the inclusion $C\leq Z(W)$. For every point $x$, we can always find two points $y$ and $z$ such that $d(y,x) = 1$, $d(x,z) = 2$ and $d(y,z) = 3$. By Corollary \ref{distance 2}, $[V_x,V_y] = [V_x,V_z] = 1$. Hence $[V_x, [V_y,V_z]] = 1$. However $C = [V_y,V_z]$, by Lemma \ref{commutatori}. Therefore $[V_x,C] = 1$. It follows that $[W, C] = 1$, since $W$ is generated by the subgroups $V_x$ for $x\in{\cal P}$. \eop  

\begin{prop}\label{final}
The group $W$ is either elementary abelian (of order $p^{56}$) or extraspecial (order $p^{1+56}$). 
\end{prop}
\pr If $Z(W) = W$, namely $W' = 1$, then $W$ is elementary abelian of order $p^{56}$, by Corollary \ref{W/W' cor}.  Suppose that $Z(W) < W$. We have $W' \leq Z(W)$ by Corollary \ref{commutatori bis}. We shall prove that $Z(W) = W'$.

We know that $W'$ defines a quotient of $\eta$ and that $\eta/W'\cong \varepsilon$ (Lemma \ref{W/W'}). By the same arguments used in the proof of Lemma \ref{W/W'} one can prove that $Z(W)\cap V_X = 1$ for every element $X\in \Gamma$ and $(Z(W)V_x)\cap(Z(W)V:y) = Z(W)$ for any two distinct points $x, y \in {\cal P}$. Therefore $Z(W)$ defines a quotient of $\eta$. Clearly, $\eta/Z(W)$ is a quotient of $\eta/W'$, since $W' \leq Z(W)$. Moreover, $\eta/Z(W)$ is homogeneous, since $\eta$ is homogeneous by assumption and $Z(W)$ is characteristic in $W$. Thus, we have obtained a homogeneous proper quotient $\eta/Z(W)$ of $\varepsilon \cong \eta/W'$. However $\varepsilon$ does not admit any proper homogeneous quotient (Corollary \ref{0}). Hence $\eta/Z(W)\cong \varepsilon \cong \eta/W'$, namely $Z(W) = W'$. 

In order to prove that $W$ is extraspecial we still must show that $W'$ is an elementary $p$-group, namely all of its non-trivial elements have order $p$. We now that $W' = [V_x,V_y]$ for a given pair of points at distance $3$ (lemma \ref{commutatori}). For $a\in V_x$ and $b\in V_y$, we shall prove that
\begin{equation}\label{potenza commutatore}
(aba^{-1}b^{-1})^k = ab^ka^{-1}b^{-k}
\end{equation}
for any integer $k\leq 0$. We shall prove this by induction on $k$. Equality (\ref{potenza commutatore}) is trivial when $0 \leq k \leq 1$. Let $k > 1$. Then 
\[(aba^{-1}b^{-1})^k = aba^{-1}b^{-1}\cdot(aba^{-1}b^{-1})^{k-1} = aba^{-1}b^{-1}\cdot ab^{k-1}a^{-1}b^{-k+1}\]
where the second equality holds by the inductive ipothesis. On the other hand,
\[aba^{-1}b^{-1}\cdot ab^{k-1}a^{-1}b^{-k+1} = aba^{-1}\cdot ab^{k-1}a^{-1}b^{-k+1}\cdot b^{-1}\]
since $b^{-1}$ commutes with $ab^{k-1}a^{-1}b^{-k+1}\in [V_x,V_y] = W' = Z(W)$. Moreover $aba^{-1}\cdot ab^{k-1}a^{-1}b^{-k+1}\cdot b^{-1} = ab^ka^{-1}b^{-k}$. Hence (\ref{potenza commutatore}) holds for the given value of $k$. By induction, (\ref{potenza commutatore}) holds for any $k$. In particular, with $k = p$ we get $(aba^{-1}b^{-1})^p = ab^pa^{-1}b^{-p}$. On the other hand, $b^p = 1$ because $|V_y| = p$. Hence  $(aba^{-1}b^{-1})^p = 1$.  Thus, $W'$ is indeed elementary.  

It is well known that the center of an extraspecial $p$-group is cyclic of order $p$. So, $|W| = p^{1+56}$, since $|W/W'| = p^{56}$.   \eop 

\bigskip

\noindent
{\bf End of the proof of Theorem \ref{ThB2}.} Let now $\Delta = E_8(p)$ and let $\varepsilon_\Delta^8$ be the unitary representation of $\Gamma$ relative to $\Delta$ and $8$. Let $a$ and $a^*$ be two opposite $8$-elements of $\Delta$. Property $(\mathrm{A})$ of Subsection \ref{sezione-epsilonA} holds for the local geometry $\Delta_8(a)$. Hence $\varepsilon_\Delta^J$ is homogeneous. Moreover, chosen a $1$-element $S$ of $\Delta$ incident to $a$, the geometry $\Delta_8(A)_S$ (see Subsection \ref{indotto sezione}) is isomorphic to $D_{6,1}(p)$ and $\varepsilon_\Delta^8$ induces on it the unipotent representation of $D_{6,1}(p)$ (compare Theorem \ref{unipotente indotto}). By claim (3) of Theorem \ref{ThA1}, the latter is isomorphic to the $12$-dimensional projective embedding of $D_{6,1}(p)$. It follows that $\varepsilon_\Delta^J$ is locally projective. (Recall that every element of $\Delta_8(a)$ is incident to at least one $1$-element.) 

Thus, $\varepsilon_\Delta^8$ satisfies both conditions $(\eta 1)$ and $(\eta 2)$. By Proposition \ref{final}, the group $U_{a^*}$ is either elementary abelian of order $p^{56}$ or extraspecial of type $p^{1+56}$. The order of $U_{a^+}$ is easy to compute. In fact, $|U_{a^*}| = p^{57}$. Therefore $U_{a^*}$ is extraspecial. 

We have $\varepsilon \cong \varepsilon_{\Delta,\mathrm{ab}}^8$ by Proposition \ref{W/W'}. The embedding $\varepsilon$ is not dominant, since it is a proper quotient of $\varepsilon_\Delta^8$. On the other hand, the codomain $U(\varepsilon)$ of the hull $\tilde{\varepsilon}$ of $\varepsilon$ cannot have order greater than $p^{57}$, by Proposition \ref{final}. It follows that $\varepsilon_\Delta^8\cong \tilde{\varepsilon}$.  Both claims of Theorem \ref{ThB2} are proved. \eop

\bigskip

\noindent
{\bf Remark 1.}  We have proved that $U_{a^*}$ is extraspecial. Actually $U_{a^*}$ is extraspecial of $+$ type, but we are not going to prove this fact here.\\

\noindent
{\bf Remark 2.} The hypothesis that $\mathbb{F} = \mathbb{F}_p$ has been used several times in our proof of Theorem \ref{ThB2}. In particular, it is used in an essential way at the end of our proof, to determine the feasible structure of $W$ (Proposition \ref{final}) and to conclude that $U_{a^*} \cong U(\varepsilon)$ by comparing orders. In the general case, with $\mathbb{F}_p$ replaced by an arbitary field, we should use quite different arguments. For instance, we could do as follows. We might firstly check that $\varepsilon \cong \varepsilon_{\Delta,\mathrm{ab}}^8$ by a direct inspection of $U_{a^*}$. Next we should prove that $\mathrm{Far}_\Delta(a^*)$ is $2$-simply connected, so that to combine this fact with Theorems \ref{1.5} and \ref{expansion} to conclude that $\varepsilon_\Delta^8$ is dominant.  

\section{Notation}\label{notation}

\subsection{Notation for Coxeter diagrams and shadow geometries}\label{notation-Coxeter}

Following a well established custom, we use positive integers as types, labelling the nodes of an irreducible spherical Coxeter diagram in the following way:

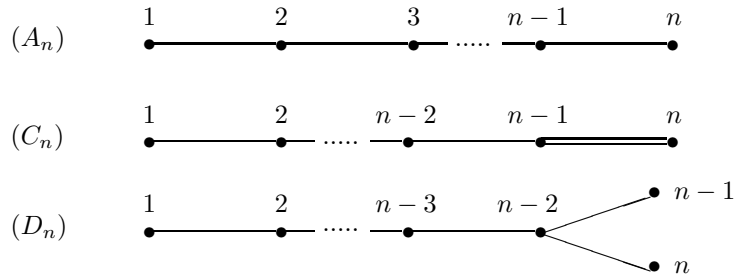
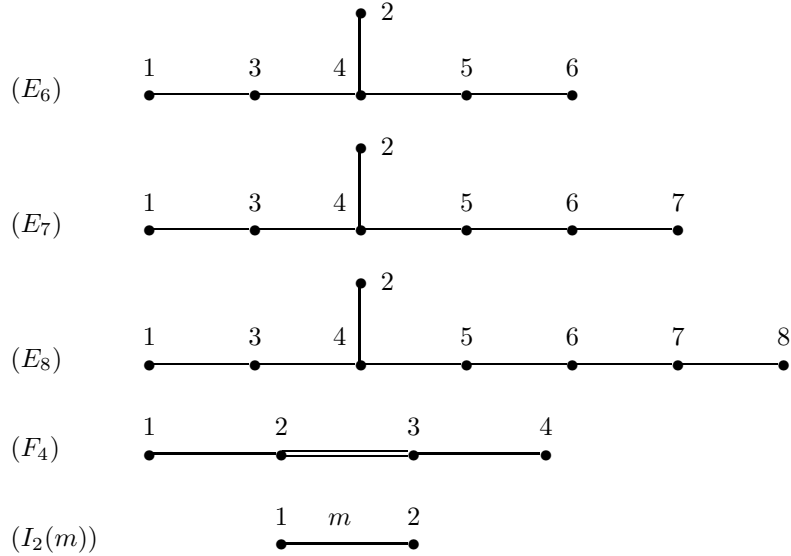
\begin{picture}(310,36)(0,0)
\put(0,10){($A_n$)}
\put(50,8){$\bullet$}
\put(53,11){\line(1,0){47}}
\put(100,8){$\bullet$}
\put(103,11){\line(1,0){47}}
\put(150,8){$\bullet$}
\put(153,11){\line(1,0){12}}
\put(168,10){$.....$}
\put(186,11){\line(1,0){12}}
\put(198,8){$\bullet$}
\put(201,11){\line(1,0){47}}
\put(248,8){$\bullet$}
\put(50,18){$1$}
\put(100,18){$2$}
\put(150,18){$3$}
\put(188,18){$n-1$}
\put(248,18){$n$}
\end{picture}

\begin{picture}(310,36)(0,0)
\put(0,10){($C_n$)}
\put(50,8){$\bullet$}
\put(53,11){\line(1,0){47}}
\put(100,8){$\bullet$}
\put(103,11){\line(1,0){12}}
\put(118,10){$.....$}
\put(136,11){\line(1,0){12}}
\put(148,8){$\bullet$}
\put(151,11){\line(1,0){47}}
\put(198,8){$\bullet$}
\put(201,10){\line(1,0){47}}
\put(201,12){\line(1,0){47}}
\put(248,8){$\bullet$}
\put(50,18){1}
\put(100,18){2}
\put(138,18){$n-2$}
\put(188,18){$n-1$}
\put(248,18){$n$}
\end{picture}

\begin{picture}(310,46)(0,0)
\put(0,23){($D_n$)}
\put(50,21){$\bullet$}
\put(53,24){\line(1,0){47}}
\put(100,21){$\bullet$}
\put(103,24){\line(1,0){12}}
\put(118,24){$.....$}
\put(136,24){\line(1,0){12}}
\put(148,21){$\bullet$}
\put(151,24){\line(1,0){47}}
\put(198,21){$\bullet$}
\put(200,23){\line(3,1){42}}
\put(200,23){\line(3,-1){42}}
\put(241,36){$\bullet$}
\put(241,8){$\bullet$}
\put(50,31){1}
\put(100,31){2}
\put(138,31){$n-3$}
\put(184,31){$n-2$}
\put(251,8){$n$}
\put(251,36){$n-1$}
\end{picture}

\begin{picture}(310,43)(0,0)
\put(0,10){($E_6$)}
\put(50,8){$\bullet$}
\put(53,11){\line(1,0){37}}
\put(90,8){$\bullet$}
\put(93,11){\line(1,0){37}}
\put(130,8){$\bullet$}
\put(133,11){\line(1,0){37}}
\put(170,8){$\bullet$}
\put(173,11){\line(1,0){37}}
\put(210,8){$\bullet$}
\put(132,10){\line(0,1){30}}
\put(130,39){$\bullet$}
\put(50,18){$1$}
\put(90,18){$3$}
\put(122,18){$4$}
\put(170,18){$5$}
\put(210,18){$6$}
\put(140,39){$2$}
\end{picture}

\begin{picture}(310,50)(0,0)
\put(0,10){($E_7$)}
\put(50,8){$\bullet$}
\put(53,11){\line(1,0){37}}
\put(90,8){$\bullet$}
\put(93,11){\line(1,0){37}}
\put(130,8){$\bullet$}
\put(133,11){\line(1,0){37}}
\put(170,8){$\bullet$}
\put(173,11){\line(1,0){37}}
\put(210,8){$\bullet$}
\put(213,11){\line(1,0){37}}
\put(250,8){$\bullet$}
\put(132,10){\line(0,1){30}}
\put(130,39){$\bullet$}
\put(50,18){$1$}
\put(90,18){$3$}
\put(122,18){$4$}
\put(170,18){$5$}
\put(210,18){$6$}
\put(250,18){$7$}
\put(140,39){$2$}
\end{picture}

\begin{picture}(310,50)(0,0)
\put(0,10){($E_8$)}
\put(50,8){$\bullet$}
\put(53,11){\line(1,0){37}}
\put(90,8){$\bullet$}
\put(93,11){\line(1,0){37}}
\put(130,8){$\bullet$}
\put(133,11){\line(1,0){37}}
\put(170,8){$\bullet$}
\put(173,11){\line(1,0){37}}
\put(210,8){$\bullet$}
\put(213,11){\line(1,0){37}}
\put(250,8){$\bullet$}
\put(253,11){\line(1,0){37}}
\put(290,8){$\bullet$}
\put(132,10){\line(0,1){30}}
\put(130,39){$\bullet$}
\put(50,18){$1$}
\put(90,18){$3$}
\put(122,18){$4$}
\put(170,18){$5$}
\put(210,18){$6$}
\put(250,18){$7$}
\put(290,18){$8$}
\put(140,39){$2$}
\end{picture}

\begin{picture}(310,33)(0,0)
\put(0,10){$(F_4)$}
\put(50,8){$\bullet$}
\put(53,11){\line(1,0){47}}
\put(100,8){$\bullet$}
\put(103,10){\line(1,0){47}}
\put(103,12){\line(1,0){47}}
\put(150,8){$\bullet$}
\put(153,11){\line(1,0){47}}
\put(200,8){$\bullet$}
\put(50,18){1}
\put(100,18){2}
\put(150,18){3}
\put(200,18){4}
\end{picture}

\begin{picture}(310,33)(0,0)
\put(0,10){$(I_2(m))$}
\put(100,8){$\bullet$}
\put(103,11){\line(1,0){47}}
\put(150,8){$\bullet$}
\put(120,18){$m$}
\put(100,18){1}
\put(150,18){2}
\end{picture}

\noindent
If $X_n$ is the symbol used to denote the Coxeter diagram of a building $\Delta$ of rank $n$ and $J$ is a subset of the set of types of $\Delta$, we say that the shadow geometry $\Delta_J$ has {\em type} $X_{n,J}$. Thus, for instance, $A_{n,2}$ is the type of the line-grassmannian of an $n$-dimensional projective space. If $\Delta$ has type $C_n$ or $D_n$ then $\Delta_1$ is the polar space associated to $\Delta$. It has type $C_{n,1}$ or $D_{n,1}$ respectively. A point-line geometry of type $C_{n,n}$ is a dual polar space. A point-line geometry of type $D_{n,n}$ ($\cong D_{n,n-1}$) is a half-spin geometry. Its local geometries are isomorphic to shadow geometries of type $A_{n-1,2}$. If $\Gamma := \Delta_J$ has type $E_{n,1}$ ($n = 6, 7$ or $8$) then its local geometries are half-spin geometries (type $D_{n-1,n-1}$). If $\Gamma$ has type $E_{n,n}$ for $n = 6, 7$  or $8$ then its local geometries are shadow geometries of type $D_{5,5}$, $E_{6,6}$ ($\cong E_{6,1}$) and $E_{7,7}$ respectively.  If $\Gamma$ has type $F_{4,1}$ or $F_{4,4}$ then its local geometries are dual polar spaces of rank $3$ (type $C_{3,3}$). 

\subsection{Notation for certain buildings}\label{notation-Dynkin} 

We recall that for $n \geq 4$ and every field (i.e. commutative division ring) $\mathbb{F}$ there is exactly one building of type $D_n$ defined over $\mathbb{F}$. We denote it by the symbol $D_n(\mathbb{F})$. Accordingly, $D_{n,J}(\mathbb{F})$ is the $J$-shadow geometry of $D_n(\mathbb{F})$. Similarly, for $6\leq n\leq 8$ we denote by $E_n(\mathbb{F})$ the unique building of type $E_n$ defined over $\mathbb{F}$ and by $E_{n,J}(\mathbb{F})$ its $J$-shadow geometry. 

A similar notation is used for buildings of type $A_n$, $C_n$ and $F_4$, except that the symbols $C_n$ and $F_4$ are now given a meaning sharper than in Subsection \ref{notation-Coxeter}. Indeed by $C_n(\mathbb{F})$ or $F_4(\mathbb{F})$ we do not mean just a building of Coxeter type $C_n$ of $F_4$ somehow related to the field $\mathbb{F}$. Now the letters $C_n$ and $F_4$ are names of Dynkin diagrams rather than Coxeter diagrams. The symbol $C_n(\mathbb{F})$ stands for the building associated to the symplectic group $\mathrm{Sp}(2n,\mathbb{F})$ while $F_4(\mathbb{F})$ is the building associated to an $\mathbb{F}$-split algebraic group of type $F_4$. Another Dynkin diagram exists besides $C_n$ with the same Coxeter shape as $C_n$, namely $B_n$: the symbol $B_n(\mathbb{F})$ denotes the building (of Coxeter type $C_n$) associated to the spin group $\mathrm{Spin}(2n+1,\mathbb{F})$.
 Note that $B_2(\mathbb{F})$ and $C_2(\mathbb{F})$ are mutually dual. It is custom to distinguish betwen them by putting the sign $>$ or $<$ on the double stroke of their Coxeter diagram $C_2$, using this trick also for cases of rank $n > 2$.          

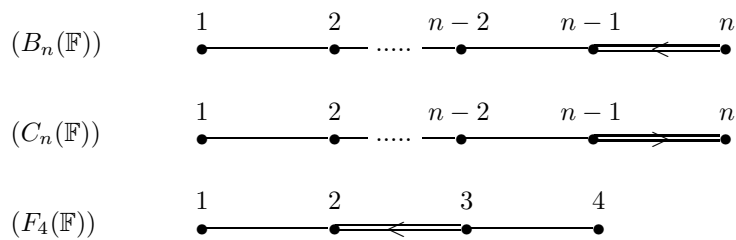
\begin{picture}(310,36)(0,0)
\put(0,10){($B_n(\mathbb{F}))$}
\put(70,8){$\bullet$}
\put(73,11){\line(1,0){47}}
\put(120,8){$\bullet$}
\put(123,11){\line(1,0){12}}
\put(138,10){$.....$}
\put(156,11){\line(1,0){12}}
\put(168,8){$\bullet$}
\put(171,11){\line(1,0){47}}
\put(218,8){$\bullet$}
\put(221,10){\line(1,0){47}}
\put(221,12){\line(1,0){47}}
\put(268,8){$\bullet$}
\put(242,8){$<$}
\put(70,18){1}
\put(120,18){2}
\put(158,18){$n-2$}
\put(208,18){$n-1$}
\put(268,18){$n$}
\end{picture}

\begin{picture}(310,33)(0,0)
\put(0,10){($C_n(\mathbb{F}))$}
\put(70,8){$\bullet$}
\put(73,11){\line(1,0){47}}
\put(120,8){$\bullet$}
\put(123,11){\line(1,0){12}}
\put(138,10){$.....$}
\put(156,11){\line(1,0){12}}
\put(168,8){$\bullet$}
\put(171,11){\line(1,0){47}}
\put(218,8){$\bullet$}
\put(221,10){\line(1,0){47}}
\put(221,12){\line(1,0){47}}
\put(268,8){$\bullet$}
\put(242,8){$>$}
\put(70,18){1}
\put(120,18){2}
\put(158,18){$n-2$}
\put(208,18){$n-1$}
\put(268,18){$n$}
\end{picture}

\begin{picture}(310,33)(0,0)
\put(0,10){$(F_4(\mathbb{F}))$}
\put(70,8){$\bullet$}
\put(73,11){\line(1,0){47}}
\put(120,8){$\bullet$}
\put(123,10){\line(1,0){47}}
\put(123,12){\line(1,0){47}}
\put(170,8){$\bullet$}
\put(173,11){\line(1,0){47}}
\put(220,8){$\bullet$}
\put(142,8){$<$}
\put(70,18){1}
\put(120,18){2}
\put(170,18){3}
\put(220,18){4}
\end{picture}

\noindent
As for buildings of twisted type, we only recall the following conventions. Let $\mathbb{F}$ be a field admitting a quadratic extension $\mathbb{K}$. We denote by ${^2}A_n(\mathbb{F})$ the building of Coxeter type $C_n$ formed by the subspaces of $V := V(n+1,\mathbb{K})$ totally isotropic for a given non-degenerate $\sigma$-hermitian form of $V$, where $\sigma$ is the unique non-trivial Galois automorphism of the extension $[\mathbb{K}:\mathbb{F}]$.  

The symbol ${^2}D_{n+1}(\mathbb{F})$ stands for the building of Coxeter type $C_n$ associated to the group $\mathrm{O}^-(2n+2,\mathbb{F})$. 

The symbol ${^2}E_6(\mathbb{F})$ denotes the building of type $F_4$ with $\{1,2,3\}$-residues isomorphic to ${^2}A_5(\mathbb{F})$ and $\{2,3,4\}$-residues isomorphic to $^{2}D_4(\mathbb{F})$.

Finally, when $\mathbb{F} = \mathbb{F}_q$ for a prime power $q$ we write $B_n(q)$ for $B_n(\mathbb{F}_q)$, $C_n(q)$ for $C_n(\mathbb{F}_q)$, and so on.

\bigskip

\noindent
Antonio Pasini,\\
Department of Information Engineering and Mathematics,\\
University of Siena, Via Roma 56, 53100 Siena, Italy.\\
antonio.pasini@unisi.it

\end{document}